%April 30, 2006
\input amstex
\documentstyle{amsppt}
\magnification=1200
\hoffset=-0.5pc
\nologo
\vsize=57.2truepc
\hsize=38.5truepc

\spaceskip=.5em plus.25em minus.20em

\define\FFF{F}
\define\ETAA{\eta}
\define\BETA{\beta}
\define\TT{T}

\define\rbieltwo{1}
\define\broetomd{2}
\define\cerezone{3}
\define\duflotwo{4}
\define\duiskolk{5}
\define\farauone{6}
\define\bhallthr{7}
\define\bhallfou{8}
\define\bhallone{9}
\define\helbotwo{10}
\define\poiscoho{11}
 \define\souriau{12}
 \define\kaehler{13}
   \define\lradq{14}
      \define\qr{15}
   \define\varna{16}
\define\bedlepro{17}
\define\adjoint{18}
\define\hurusch{19}
\define\kirilboo{20}
\define\kirilthr{21}
\define\lassaone{22}
\define\lassatwo{23}
\define\lempszoe{24}
\define\liuwanhu{25}
\define\enelson{26}
\define\sniabook{27}
\define\tspritwo{28}
\define\stenzone{29}
\define\szoekone{30}
\define\taylothr{31}
\define\woodhous{32}

\topmatter
\title Kirillov's character formula, the holomorphic Peter-Weyl theorem, and
 the Blattner-Kostant-Sternberg pairing
\endtitle
\author Johannes Huebschmann \footnote"*"{Support by the
German Research Council (Deutsche Forschungsgemeinschaft) in the
framework of a Mercator visiting professorship is gratefully
acknowledged. \hfill \hfill}
\endauthor

\date{April 30, 2007} \enddate
\abstract{Let $K$ be a compact Lie group, endowed with a
bi-invariant Riemannian metric, which we denote by $\kappa$. The
complexification $K^{\Bbb C}$ of $K$ inherits a K\"ahler structure
having twice the kinetic energy of the metric as its potential;
let $\varepsilon$ denote the symplectic volume form. Left and
right translation turn the Hilbert space $\Cal H L^2(K^{\Bbb
C},\roman e^{-\kappa /t}\eta\varepsilon)$ of square-integrable
holomorphic functions on $K^{\Bbb C}$ relative to a suitable
measure
 written as
$e^{-\kappa /t}\eta\varepsilon$ into a unitary $(K \times
K)$-representation; here   $\eta$ is an additional term coming
from the metaplectic correction, and $t>0$ is a real parameter. In
the physical interpretation, this parameter amounts to {\it
Planck\/}'s constant $\hbar$.

\vskip1ex

We establish the statement of the Peter-Weyl theorem for the
Hilbert space $\Cal H L^2(K^{\Bbb C},\roman e^{-\kappa
/t}\eta\varepsilon)$ to the effect that

 (i) $\Cal H
L^2(K^{\Bbb C}, \roman e^{-\kappa /t} \eta \varepsilon)$ contains
the vector space of representative functions on $K^{\Bbb C}$ as a
dense subspace and that

(ii) the assignment to a holomorphic function of its Fourier
coefficients yields an isomorphism of Hilbert algebras from the
convolution algebra $\Cal H L^2(K^{\Bbb C}, \roman e^{-\kappa /t}
\eta \varepsilon)$ onto an algebra of the kind $\widehat\oplus
\roman{End}(V)$. Here $V$ ranges over the irreducible rational
representations of $K^{\Bbb C}$ and $\widehat\oplus$ refers to a
suitable completion of the direct sum algebra $\oplus
\roman{End}(V)$. \vskip1ex

Consequences are:

(i) the existence of a
uniquely determined unitary isomorphism between $L^2(K,dx)$ (where
$dx$ refers to Haar measure on $K$) and the Hilbert space $\Cal H
L^2(K^{\Bbb C},\roman e^{-\kappa /t}\eta\varepsilon)$, and

(ii) a proof that this isomorphism coincides with the
Blattner-Kostant-Sternberg pairing map from  $L^2(K,dx)$  to $\Cal
H L^2(K^{\Bbb C},\roman e^{-\kappa /t}\eta\varepsilon)$,
multiplied by $(4\pi t)^{-\dim (K)/4}$. 

\vskip1ex 
Among our crucial tools is Kirillov's
character formula. Our methods are geometric,
rely on the orbit method, and are independent of
heat kernel harmonic analysis, which is used 
by B. C. Hall to obtain many of these
results [J. of Funct. Anal.  122 (1994), 103--151],
[Comm. in Math. Physics 226 (2002),  233--268]. }
\endabstract

\address{\smallskip
\noindent
USTL, UFR de Math\'ematiques, CNRS-UMR 8524
\newline\noindent
59655 VILLENEUVE  d'ASCQ C\'edex, France
\newline\noindent
Johannes.Huebschmann\@math.univ-lille1.fr}
\endaddress
\endkeywords

\endtopmatter
\document
\leftheadtext{Johannes Huebschmann}
\rightheadtext{Holomorphic Peter-Weyl theorem}

\medskip\noindent {\bf Introduction}
\smallskip\noindent
Let $K$  be a compact Lie group and let $K^{\Bbb C}$ be its
complexification. Given a finite dimensional rational
representation $V$ of $K^{\Bbb C}$, the familiar assignment to
$$\varphi \otimes w \in V^* \otimes V\cong \roman{End}(V)$$ of the
representative function $\Phi_{\varphi,w}$ given by
$\Phi_{\varphi,w}(q) = \varphi(qw)$ ($q \in K^{\Bbb C}$) yields an
embedding of $\roman{End}(V)$ into the Hopf algebra $\Bbb
C[K^{\Bbb C}]$ of representative functions on $K^{\Bbb C}$, the
diagonal map of $\Bbb C[K^{\Bbb C}]$ being induced by the group
multiplication. As $V$ ranges over the finite dimensional
irreducible representations of $K^{\Bbb C}$, these maps assemble
to an isomorphism
$$
\oplus \roman{End}(V) \longrightarrow \Bbb C[K^{\Bbb C}] \tag0.1
$$
of $(K^{\Bbb C}\times K^{\Bbb C})$-representations. Restriction to
$K$ induces an isomorphism from the Hopf algebra $\Bbb C[K^{\Bbb
C}]$ onto the Hopf algebra $R(K)$ of complex representative
functions of $K$ and induces an isomorphism
$$
\oplus \roman{End}(V) \longrightarrow R(K) \tag0.2
$$
of $(K \times K)$-representations. Integration over $K$ relative
to Haar measure induces inner products on the left- and right-hand
sides of (0.2), the completion of the right-hand side is the
Hilbert space $L^2(K,dx)$ relative to Haar measure $dx$ which, as
a unitary $K$-representation, is the {\it regular\/}
representation, and the completion
$$
\widehat\oplus \roman{End}(V) \longrightarrow L^2(K,dx) \tag0.3
$$
of the isomorphism (0.2) yields the familiar {\it Peter-Weyl\/}
theorem for the compact group $K$. Here and below the notation
$\widehat\oplus$ refers to a Hilbert space direct sum involving
infinitely many summands. Moreover, integration over $K$ relative
to Haar measure, suitably normalized, induces a convolution
product on $R(K)$ and on $L^2(K,dx)$.  Relative to this
convolution product, when $\oplus \roman{End}(V)$ is endowed with
the obvious algebra structure having the $\roman{End}(V)$'s as
minimal two-sided ideals, (0.2) is an isomorphism of algebras
(without 1 unless $K$ is a finite group), and (0.3) is an
isomorphism of Hilbert algebras and hence yields in particular the
decomposition of the convolution algebra $L^2(K,dx)$ into minimal
two-sided topological ideals. Furthermore, a choice of
bi-invariant Riemannian metric on $K$ determines a
Laplace-Beltrami operator on $K$ which admits a unique extension
to an (unbounded) self-adjoint operator on $L^2(K,dx)$, and the
decomposition on the left-hand side of (0.3) is precisely the
standard refinement of the spectral decomposition of this
operator. This extension of the Laplace-Beltrami operator on $K$
to an operator on $L^2(K,dx)$, multiplied by $-1/2$, is the
(quantum mechanical) energy operator on $L^2(K,dx)$ associated
with the metric; indeed, vertical dequantization of this operator
or, equivalently, the operation of passing to the symbol, yields
the energy associated with the Riemannian metric on $K$, viewed as
a function on the total space $\roman T^*K$ of the cotangent
bundle of $K$.

Since $K^{\Bbb C}$ is reductive, the coordinate ring of $K^{\Bbb
C}$ coincides with the algebra of representative functions on
$K^{\Bbb C}$, and hence the Hilbert space $L^2(K,dx)$ contains the
complex vector space underlying the coordinate ring of $K^{\Bbb
C}$. However, a more natural Hilbert space containing the
coordinate ring would be a Hilbert space of holomorphic functions
on $K^{\Bbb C}$. This raises the question whether $K^{\Bbb C}$
carries a suitable measure such that, after completion relative to
this measure, (0.1) yields an isomorphism of Hilbert spaces
between a Hilbert space of the kind $\widehat\oplus
\roman{End}(V)$ and a Hilbert space of holomorphic functions on
$K^{\Bbb C}$, and what the spectral decomposition of the energy
operator might correspond to in terms of $K^{\Bbb C}$. Actually,
we were led to these questions by the observation that K\"ahler
quantization is, perhaps, better suited to explore quantization in
the presense of singularities than ordinary Schr\"odinger
quantization. We shall comment on this motivation below.

The desired measure on $K^{\Bbb C}$ is provided for by the measure
coming from half-form quantization on $K^{\Bbb C}$, though the
construction of the measure itself is independent of the program
of geometric quantization: A choice of bi-invariant Riemannian
metric on $K$ and the polar decomposition map of $K^{\Bbb C}$
determine a diffeomorphism between $\roman T^*K$ and $K^{\Bbb C}$
and, via this diffeomorphism, $\roman T^*K$ and $K^{\Bbb C}$ both
become K\"ahler manifolds where the requisite symplectic structure
is the cotangent bundle structure on $\roman T^*K$. Since $K^{\Bbb
C}$ is parallelizable, it admits a metaplectic structure for
trivial reasons, and the bundle of holomorphic half-forms on
$K^{\Bbb C}$ furnishes a measure on $K^{\Bbb C}$ having the
desired properties. This measure can be written in the form
$\roman e^{-\kappa /t} \eta \varepsilon$; here $\kappa$ is the
metric, written as a function on $\roman T^*K$, $t>0$ is a real
parameter which, in the physical interpretation, amounts to {\it
Planck\/}'s constant $\hbar$,  $\varepsilon$ is the Liouville
volume measure, and  $\eta$ is a suitable function coming from the
metaplectic correction and closely related with the familiar
function coming into play in {\it Kirillov\/}'s character formula
and commonly written as $j$ \cite\duflotwo, \cite\kirilthr.

We shall establish a holomorphic version of the {\it Peter-Weyl\/}
theorem to the effect that the following hold: (i) {\sl The
Hilbert space $\Cal H L^2(K^{\Bbb C}, \roman e^{-\kappa /t} \eta
\varepsilon)$ of holomorphic functions that are square integrable
relative to the measure $\roman e^{-\kappa /t} \eta \varepsilon$
contains the vector space $\Bbb C[K^{\Bbb C}]$ of representative
functions on $ K^{\Bbb C}$ as a dense subspace in such a way that
the decomposition\/} (0.1) {\sl induces the decomposition of $\Cal
H L^2(K^{\Bbb C}, \roman e^{-\kappa /t} \eta \varepsilon)$, viewed
as a unitary $(K\times K)$-representation, into its isotypical
summands\/}; and  (ii) {\sl the assignment to a holomorphic
function of its Fourier coefficients yields an isomorphism of
Hilbert algebras from $\Cal H L^2(K^{\Bbb C}, \roman e^{-\kappa
/t} \eta \varepsilon)$, made into an algebra via the convolution
product, onto an algebra of the kind $\widehat\oplus
\roman{End}(V)$ where $V$ ranges over the irreducible rational
representations of $K^{\Bbb C}$ and where $\widehat\oplus$ refers
to a suitable completion of the direct sum algebra $\oplus
\roman{End}(V)$\/}; a precise statement is given as Theorem 1.14
below. We then refer to the resulting decomposition as the {\it
holomorphic Peter-Weyl decomposition\/} of the Hilbert space $\Cal
H L^2(K^{\Bbb C}, \roman e^{-\kappa /t} \eta \varepsilon)$. The
algebraic decomposition (0.1) itself into isotypical $(K^{\Bbb
C}\times K^{\Bbb C})$-summands is commonly interpreted as an {\it
algebraic Peter-Weyl\/} theorem.

Consequences of the holomorphic Peter-Weyl theorem are
 the existence of a uniquely determined unitary
$(K\times K)$-equivariant isomorphism between $L^2(K,dx)$ (where
$dx$ refers to Haar measure on $K$) and the Hilbert space $\Cal H
L^2(K^{\Bbb C},\roman e^{-\kappa /t}\eta\varepsilon)$, given in
Theorem 5.3 below, and a {\it holomorphic Plancherel\/} theorem,
given as Corollary 5.4 below. In Section 6 we shall then show that
the isomorphism between the two Hilbert spaces coincides with the
Blattner-Kostant-Sternberg pairing map from the Hilbert space
$L^2(K,dx)$ to the Hilbert space $\Cal H L^2(K^{\Bbb C},\roman
e^{-\kappa /t}\eta\varepsilon)$, multiplied by $(4\pi t)^{-\dim
(K)/4}$. However the abstract unitary isomorphism between the two
Hilbert spaces is independent of the Blattner-Kostant-Sternberg
pairing. The identification of the two Hilbert spaces implies, in
particular, that the spectral decomposition of the energy operator
on $\Cal H L^2(K^{\Bbb C},\roman e^{-\kappa /t}\eta\varepsilon)$
associated with the metric refines to the holomorphic Peter-Weyl
decomposition of this Hilbert space in the usual manner and thus
yields the decomposition of $\Cal H L^2(K^{\Bbb C},\roman
e^{-\kappa /t}\eta\varepsilon)$ into irreducible isotypical $(K
\times K)$-representations; this will be explained in Section 7.

A crucial step towards the holomorphic Peter-Weyl theorem consists
in proving that the representative functions on $K^{\Bbb C}$ are
square integrable relative to the measure $\roman e^{-\kappa /t}
\eta \varepsilon$ on $K^{\Bbb C}$ and, furthermore, in actually
calculating their square integrals. We do these calculations by
means of {\it Kirillov\/}'s character formula; we could as well
have taken {\it Weyl\/}'s character formula, but then the
calculations would be somewhat more involved. Our argument
establishing the completeness of the representative functions,
given in Section 4 below, is geometric and is guided by the
principle that quantization commutes with reduction.

An obvious question emerges here: How is the Hilbert space $\Cal H
L^2(K^{\Bbb C},\roman e^{-\kappa /t}\varepsilon)$ of holomorphic
functions that are square integrable relative to the measure
$\roman e^{-\kappa /t}\varepsilon$ related to the other Hilbert
spaces? In a final section we shall show that, indeed, as a
unitary $(K\times K)$-representation, this Hilbert space is
unitarily equivalent to the Hilbert space $\Cal H L^2(K^{\Bbb
C},\roman e^{-\kappa /t}\eta\varepsilon)$ in an obvious manner.

We now relate the present paper with what we know to be already in
the literature. Since the measure $\roman e^{-\kappa
/t}\eta\varepsilon$ involves a Gaussian constituent, the
square-integrability of the representative functions relative to
the corresponding measure can also be established directly, and
Lemma 10 in \cite\bhallthr, combined with the observation just
made, entails that the representative functions are dense in $\Cal
H L^2(K^{\Bbb C},\roman e^{-\kappa /t}\eta\varepsilon)$. The fact
that the BKS-pairing map between the two Hilbert spaces,
multiplied by a suitable constant, is a unitary isomorphism, has
been established by B. Hall \cite\bhallone. In that paper, the
pairing map is shown to coincide, up to multiplication by  a
constant, with a version of the Segal-Bargmann coherent state
transform developed, in turn, over Lie groups admitting a
bi-invariant Riemannian metric,  in a sequence of preceding papers
\cite\bhallthr --\cite\bhallone. The main technique in those
papers is heat kernel harmonic analysis and, in fact, in
\cite\bhallone, Hall derives the unitarity of the pairing map by
identifying the measure on $K^{\Bbb C}$ coming from the half-form
bundle with an appropriate heat kernel measure which, in turn, he
has shown in the preceding papers to furnish a unitary transform.
This approach in terms of the Segal-Bargmann transform, combined
with the ordinary Peter-Weyl theorem, also entails the statement
of the holomorphic Peter-Weyl theorem. A version of the
holomorphic Plancherel Theorem may be found in \cite\lassaone\ as
well as in Lemmata 9 and 10 of \cite\bhallthr. In \cite\bhallthr\
(Section 8), the completeness of the representative functions is
established by analytical considerations. A number of results in
Section 10 of \cite\bhallthr\ are actually independent of heat
kernel methods in the sense that they are valid for more general
measures than that coming from heat kernel analysis and, when
these results are applied to the heat kernel measure, evaluation
of coefficients is possible in terms of the eigenvalues of the
Laplacian. The uniform convergence, on compact sets of a group of
the kind $K^{\Bbb C}$, of what we refer to as the holomorphic
Fourier series may be found already in Proposition 12 of
\cite\cerezone. See also Remark 5.5 below.

Our methods are direct and independent of heat kernels and of the
Laplacian, involve little analysis, if any, and imply, in
particular, that the unitary isomorphism between the two Hilbert
spaces is independent of heat kernels; see Remarks 6.8 and 7.4
below. Thus, our approach answers a question raised in \cite
\bhallone, see Remark 5 in (2.5) of \cite \bhallone; it also paves
the way towards exploring the unitarity issue of the BKS-pairing
map over homogeneous spaces, to which we will come back elsewhere.
The Segal-Bargmann transform on a symmetric space of compact type
has been studied in \cite\stenzone\ but the question how this
transform is related with the corresponding BKS-pairing has not
been investigated in that paper. We plan to show at another
occasion that the description of the pairing maps in terms of heat
kernels is a direct consequence of our method. For comparison of
the two approaches with the BKS-pairing, see the identity (7.2)
below. It is, perhaps, also worthwhile pointing out that the space
$\Cal H L^2(K^{\Bbb C},\roman e^{-\kappa /t}\eta\varepsilon)$ is a
{\it weighted Bergman space\/} but we shall not use the theory of
general Bergman spaces.

This paper was written during a stay at the Institute for
Theoretical Physics at the University of Leipzig. This stay was
made possible by the German Research Council (Deutsche
Forschungsgemeinschaft) in the framework of a Mercator visiting
professorship, and I wish to express my gratitude to this
organization. It is a pleasure to acknowledge the stimulus of
conversation with G. Rudolph and M. Schmidt at Leipzig. The paper
is part of a research program aimed at exploring quantization on
classical phase spaces with singularities \cite\poiscoho
--\cite\bedlepro, in particular on classical lattice gauge theory
phase spaces. Details for the special case of a single spatial
plaquette where $K=\roman{SU}(2)$ are worked out in \cite\hurusch.
The precise information needed for this research program is the
equivalence of the two Hilbert spaces spelled out in Theorem 5.3
below.

I am indebted to R. Sz\"oke, B.~Hall and M. Lassalle for
discussion and for having provided important information which
helped placing the paper properly in the literature.

\beginsection 1. The Peter-Weyl decomposition of the half-form Hilbert
space

Let $K$ be a compact Lie group and $K^{\Bbb C}$ its
complexification, and let $\frak k$ and $\frak k^{\Bbb C}$ be the
Lie algebras of $K$ and $K^{\Bbb C}$, respectively. Choose an
invariant inner product $\cdot \, \colon \frak k \otimes \frak k
@>>> \Bbb R$ on $\frak k$, and endow $K$ with the corresponding
bi-invariant Riemannian metric. Using the metric, we identify
$\frak k$ with its dual $\frak k^*$ and the total space $\roman T
K$ of the tangent bundle with the total space  $\roman T^* K$ of
the cotangent bundle, and we will denote by $|\,\cdot\,|$ the
resulting norms on $\frak k$ and on $\frak k^*$.

Consider the {\it polar decomposition\/} map
$$
K\times \frak k \longrightarrow K^{\Bbb C},\ (x,Y) \mapsto
x\cdot\roman{exp}(iY), \ (x,Y) \in K \times \frak k. \tag1.1
$$
The composite of the inverse of left trivialization with (1.1)
identifies $\roman T^* K$ with $K^{\Bbb C}$ in a $(K\times
K)$-equivariant fashion. Then the induced complex structure on
$\roman T^* K$ combines with the symplectic structure to a
$K$-bi-invariant (positive) K\"ahler structure. Indeed, the {\it
real analytic\/} function
$$
\kappa \colon K^{\Bbb C} @>>> \Bbb R, \quad
\kappa(x\cdot\roman{exp}(iY)) =  |Y|^2,\quad (x,Y) \in K \times
\frak k, \tag1.2
$$
on $K^{\Bbb C}$ which is twice the {\it kinetic energy\/}
associated with the Riemannian metric, is a (globally defined)
$K$-bi-invariant {\it K\"ahler potential\/}; in other words, the
function $\kappa$ is strictly plurisubharmonic and (the negative
of the imaginary part of) its Levi form yields (what corresponds
to) the cotangent bundle symplectic structure, that is, the
tautological cotangent bundle symplectic structure on $\roman
T^*K\cong K^{\Bbb C}$ is given by $i\,\partial \overline\partial
\kappa$. An explicit calculation which establishes this fact may
be found in \cite\bhallone.  For related questions see
\cite\lempszoe, \cite\szoekone.

We now introduce an additional real parameter $t>0$; in the
physical interpretation, this parameter amounts to {\it
Planck\/}'s constant $\hbar$. For ease of comparison with the heat
kernel measure, cf. the identity (7.3) below, we prefer the
notation $t$ rather than $\hbar$. Half-form K\"ahler quantization,
cf. e.~g. \cite\woodhous\ (chap. 10), applied to $K^{\Bbb C}$
relative to the tautological cotangent bundle symplectic structure
on $K^{\Bbb C}$, multiplied by $1/t$, is accomplished by means of
a certain Hilbert space of holomorphic functions on $K^{\Bbb C}$
which we now recall, for ease of exposition; see \cite\bhallone\
for details. For the sake of brevity, we do not spell out the
half-forms explicitly.

Let $\varepsilon$ be the symplectic (or Liouville) volume form on
$\roman T^*K\cong K^{\Bbb C}$; this form induces the Liouville
volume measure, and we will refer to $\varepsilon$ as {\it
Liouville (volume) measure\/} as well. Further, let $dx$ denote
the volume form on $K$ yielding {\it Haar\/} measure, {\it
normalized\/} so that it coincides with the {\it Riemannian volume
measure\/} on $K$,
 and let $dY$ be the volume form inducing
{\it Lebesgue\/} measure on $\frak k$, {\it normalized\/} by the
{\it inner product\/} on $\frak k$. In terms of the polar
decomposition (1.1), we then have the identity $\varepsilon=dxdY$.
We prefer {\it not\/} to normalize the inner product on $\frak k$
since this inner product yields the kinetic energy.

Define the function $ \ETAA \colon K^{\Bbb C} \longrightarrow \Bbb
R$ by
$$
\ETAA(x,Y) = \left(\roman{det}\left(\frac{\sin(\roman{ad}(Y))}
{\roman{ad}(Y)}\right)\right)^{\frac 12}, \ x \in K, \,Y \in \frak
k; \tag1.3
$$
this yields a non-negative real analytic function on $K^{\Bbb C}$
which depends only on the variable $Y \in \frak k$ and, for $x\in
K$ and $Y\in \frak k$, we will also write $\ETAA(Y)$ instead of
$\ETAA(x,Y)$. The function $\ETAA^2$ is the density of Haar
measure relative to the Liouville volume measure on $K^{\Bbb C}$,
cf. \cite\bhallfou\ (Lemma 5). Both measures are $K$-bi-invariant;
in particular, as a function on $\frak k$, $\ETAA$ is
$\roman{Ad}(K)$-invariant. For later reference we point out that,
with the notation
$$
j(Y) =\roman{det}\left(\frac{\roman{sinh}(\roman{ad}(Y/2
))}{\roman{ad}(Y/2)}\right)^{\frac 12},\ Y\in \frak g, \tag1.4
$$
where $\frak g$ is a general Lie algebra, $ j(iY) = \eta(Y/2)$
($Y\in \frak k$). The notation $j$ is due to \cite\duflotwo\ and
\cite\kirilthr\ ((2.3.6) p.~459). We also note that a variant of
the function $\eta$ is known in the literature as the {\it van
Vleck-Morette\/} determinant. On the space of holomorphic
functions on $K^{\Bbb C}$, we will denote by $\langle \,\cdot\, ,
\, \cdot \,\rangle_{t,K^{\Bbb C}}$ the normalized inner product
given by
$$
\langle \Phi, \Psi\rangle_{t,K^{\Bbb C}} = \frac
1{\roman{vol}(K)}\int_{K^{\Bbb C}} \overline \Phi \Psi \roman
e^{-\kappa/t}\eta\varepsilon, \tag1.5
$$
and we denote by $\Cal H L^2(K^{\Bbb C},\roman e^{-\kappa/t} \ETAA
\varepsilon)$ the resulting Hilbert space of holomorphic functions
that are square integrable with respect to the measure $\roman
e^{-\kappa/t} \ETAA \varepsilon$. This
 Hilbert space is intrinsically a Hilbert space of
holomorphic half-forms on $K^{\Bbb C}$, cf. \cite\bhallone,
\cite\sniabook, \cite\woodhous. It is, furthermore, a unitary
$(K\times K)$-representation in an obvious fashion.

Given two holomorphic functions $\Phi$ and $\Psi$ on $K^{\Bbb C}$,
we define their convolution $\Phi *\Psi$ by
$$
(\Phi * \Psi)(q)  = \frac 1{\roman{vol}(K)} \int_{K}\Phi(x)
\Psi(x^{-1}q) dh_K(x),\ q \in K^{\Bbb C};\tag1.6
$$
since $K$ is compact, the convolution $\Phi *\Psi$ is a
holomorphic function on $K^{\Bbb C}$, indeed the unique extension
to a holomorphic function on $K^{\Bbb C}$ of the convolution
$(\Phi|_K)*(\Psi|_K)$ of the restrictions to $K$. Since
restriction to $K$ yields an isomorphism from $\Bbb C[K^{\Bbb C}]$
onto the space $R(K)$ of representative functions on $K$ and since
the operation of convolution turns $R(K)$ into an algebra, indeed,
a topological algebra relative to the inner product determined by
Haar measure on $K$, the operation of convolution turns the vector
space $\Bbb C[K^{\Bbb C}]$ of representative functions on $K^{\Bbb
C}$ into an algebra. We will refer to the vector space $\Bbb
C[K^{\Bbb C}]$ of representative functions on $K^{\Bbb C}$, turned
into an algebra via the convolution product, as the {\it
convolution algebra of representative functions\/} on $K^{\Bbb
C}$.

Let $T$ be a maximal torus in $K$, $\frak t$ its Lie algebra,
$T^{\Bbb C}\subseteq K^{\Bbb C}$  the complexification of $T$,
$\frak t^{\Bbb C}$ the complexification of $\frak t$, and let $W$
denote the {\it Weyl\/} group. Choose a dominant Weyl chamber
$C^+$, and let $R^+$ be the corresponding system of positive {\it
real\/} roots. Here and below the convention is that, given $Z \in
\frak t$ and an element $A$ of the root space $\frak k_{\alpha}$
associated with the root $\alpha$, the bracket $[Z,A]$ is given by
$[Z,A]= i\,\alpha(Z)A$ so that, in particular, $\alpha$ is a real
valued linear form on $\frak t$. In \cite\broetomd\ (V.1.3 on
p.~185) these $\alpha$'s are called {\it infinitesimal\/} roots.
Relative to the chosen dominant Weyl chamber, let $\rho = \frac 12
\sum_{\alpha \in R^+} \alpha$, so that $2\rho$ is the sum of the
positive roots.

We will denote by $\widehat {K^{\Bbb C}}$ the set of isomorphism
classes of irreducible rational representations of $K^{\Bbb C}$.
As usual, we identify $K^{\Bbb C}$ with the space of highest
weights relative to the chosen dominant Weyl chamber. For a
highest weight $\lambda$, we denote by $\TT_{\lambda}\colon
{K^{\Bbb C}} \to \roman{End}(V_{\lambda})$ a representation in the
class of $\lambda$ and  by $d_{\lambda}$ the dimension of
$V_{\lambda}$.

Let $\lambda$ be a highest weight. For $\psi \in V_{\lambda}^*$
and $w\in V_{\lambda}$, the function $\Phi_{\psi,w}$ given by
$$
\Phi_{\psi,w}(q) = \psi(qw),\ q \in K^{\Bbb C}, \tag 1.7
$$
is a representative function on $K^{\Bbb C}$, and the assignment
to $\psi \otimes w \in V_{\lambda}^* \otimes V_{\lambda}$ of the
representative function $\Phi_{\psi,w}$ yields a morphism
$$
\iota_{\lambda}\colon V_{\lambda}^* \otimes V_{\lambda}
\longrightarrow \Bbb C [K^{\Bbb C}]\tag1.8
$$
of $(K^{\Bbb C}\times K^{\Bbb C})$-representations, necessarily
injective  since $V_{\lambda}^* \otimes V_{\lambda}$ is an
irreducible $(K^{\Bbb C}\times K^{\Bbb C})$-representation. We
will write
$$
 V^*_{\lambda} \odot V_{\lambda}
= \iota_{\lambda}(V^*_{\lambda} \otimes V_{\lambda}) \subset \Bbb
C [K^{\Bbb C}]. \tag1.9
$$
Given an $L^2$-function $f$ on $K$ and the irreducible
representation $\TT_{\lambda}\colon K \to
\roman{End}(W_{\lambda})$ of $K$ associated with $\lambda$,
following one of the possible conventions, we define the {\it
Fourier\/} coefficient $\widehat f_{\lambda}\in
\roman{End}(W_{\lambda})$ of $f$ relative to $\lambda$ by
$$
\widehat f_{\lambda} = \frac 1{\roman{vol}(K)} \int_K
f(x)\TT_{\lambda}(x^{-1}) dx. \tag1.10
$$
Given a holomorphic function $\Phi$ on $K^{\Bbb C}$ and the
irreducible rational representation $\TT_{\lambda}\colon K^{\Bbb
C} \to \roman{End}(V_{\lambda})$ of $K^{\Bbb C}$ associated with
$\lambda$, we define the {\it Fourier\/} coefficient $\widehat
\Phi_{\lambda}\in \roman{End}(V_{\lambda})$ of $\Phi$ relative to
$\lambda$ to be the Fourier coefficient of the restriction of
$\Phi$ to $K$. The notational distinction between $V_{\lambda}$
and $W_{\lambda}$ will be justified in Section 5 below.

Let
$$
C_{t,\lambda} =(t\pi)^{\dim(K)/2}\roman e^{t|\lambda+\rho|^2}.
\tag1.11
$$
The precise significance of the real constant $C_{t,\lambda}$ will
be explained in Lemma 3.3 below. On $\roman{End}(V_{\lambda})$, we
take the standard inner product $\langle \,\cdot\, , \,\cdot\,
\rangle_{\lambda}$ given by
$$
\langle A,B\rangle_{\lambda} = \roman{tr}(A^* B),\ A,B \in
\roman{End}(V_{\lambda}), \tag1.12
$$
the adjoint $A^*$ of $A$ being computed as usual with respect to a
$K$-invariant inner product on $V_{\lambda}$. We endow ${\oplus}
_{\lambda \in \widehat{K^{\Bbb C}}}\roman{End}(V_\lambda)$ with
the  inner product which, on the summand $\roman{End}(V_\lambda)$,
is given by
$$
\frac {d_{\lambda}}{C_{t,\lambda}}\langle \,\cdot\, , \,\cdot\,
\rangle_{\lambda}; \tag1.13
$$
then $\widehat{\oplus} _{\lambda \in \widehat{K^{\Bbb
C}}}\roman{End}(V_\lambda)$ refers to the completion relative to
this inner product. Thus, up to a constant, the resulting norm on
each $\roman{End}(V_\lambda)$ is the familiar {\it
Hilbert-Schmidt\/} norm.

\proclaim{Theorem 1.14} {\rm [Holomorphic Peter-Weyl theorem]}
\newline\noindent {\rm (i)} The Hilbert space $\Cal HL^2(K^{\Bbb
C},\roman e^{-\kappa/t}\ETAA \varepsilon)$ contains the vector
space $\Bbb C [K^{\Bbb C}]$ of representative functions on
$K^{\Bbb C}$ as a dense subspace and, as a unitary $(K\times
K)$-representation, $\Cal HL^2(K^{\Bbb C},\roman
e^{-\kappa/t}\ETAA \varepsilon)$
 decomposes as
the direct sum
$$
 \Cal HL^2(K^{\Bbb
C},\roman e^{-\kappa/t}\ETAA \varepsilon)= \widehat\oplus_{\lambda
\in \widehat{K^{\Bbb C}}} V^*_{\lambda} \odot V_{\lambda}
 \tag1.14.1
$$
into $(K\times K)$-isotypical summands.
\newline\noindent
{\rm (ii)} The operation of convolution induces a convolution
product $*$ on $\Cal HL^2(K^{\Bbb C},\roman e^{-\kappa/t}\ETAA
\varepsilon)$ and, relative to this convolution product, as
$\lambda$ ranges over the irreducible rational representations of
$K^{\Bbb C}$, the assignment to a holomorphic function $\Phi$ on
$K^{\Bbb C}$ of its Fourier coefficients $\widehat
\Phi_{\lambda}\in \roman{End}(V_\lambda)$ yields an isomorphism
$$
\Cal HL^2(K^{\Bbb C},\roman e^{-\kappa/t}\ETAA \varepsilon)
\longrightarrow \widehat{\oplus} _{\lambda \in \widehat{K^{\Bbb
C}}}\roman{End}(V_\lambda) \tag1.14.2
$$
of Hilbert algebras, where each summand $\roman{End}(V_{\lambda})$
is endowed with its obvious algebra structure.
\endproclaim

The decomposition (1.14.1) of $\Cal H L^2(K^{\Bbb C},\roman
e^{-\kappa/t} \ETAA \varepsilon)$  is the {\it Peter-Weyl\/}
decomposition of this Hilbert space alluded to earlier.

\medskip\noindent
{\bf 2. The convolution algebra of representative functions}
\smallskip\noindent
For ease of exposition we recall the familiar decomposition into
minimal two-sided ideals of the convolution algebra of
representative functions on $K^{\Bbb C}$.

The operation
$$
\roman L_x\colon \Bbb C[K^{\Bbb C}] \longrightarrow \Bbb C[K^{\Bbb
C}], \ (\roman L_x (\Phi))(q) = \Phi(x^{-1}q), \ x,q \in K^{\Bbb
C}, \ \Phi \in \Bbb C[K^{\Bbb C}],
$$
of left translation on $K^{\Bbb C}$ and  the operation
$$
\roman R_y\colon \Bbb C[K^{\Bbb C}] \longrightarrow \Bbb C[K^{\Bbb
C}], \ (\roman R_y (\Phi))(q) = \Phi(qy), \ y,q \in K^{\Bbb C}, \
\Phi \in \Bbb C[K^{\Bbb C}],
$$
of right translation on $K^{\Bbb C}$ are well known to turn $\Bbb
C[K^{\Bbb C}]$ into an algebraic $(K^{\Bbb C}\times K^{\Bbb
C})$-representation in such a way that the operations of left and
right translation commute. Furthermore, the assignment to the two
representative functions $f$ and $h$ on $K^{\Bbb C}$ of $\langle
f,h\rangle = f*h(e)$ yields a complex symmetric $K^{\Bbb
C}$-invariant bilinear form
$$
\langle \,\cdot \, , \, \cdot \, \rangle \colon \Bbb C[K^{\Bbb C}]
\otimes \Bbb C[K^{\Bbb C}] \longrightarrow \Bbb C \tag2.1
$$
on $\Bbb C[K^{\Bbb C}]$, cf. e.~g. \cite\tspritwo.

Let $\lambda$ be a highest weight. We endow $ V^*_{\lambda} \odot
V_{\lambda} \cong V^*_{\lambda} \otimes V_{\lambda}$ with the
obvious complex symmetric bilinear form coming from the evaluation
mapping. By construction, this form coincides with the restriction
of the complex symmetric bilinear form {\rm (2.1)} to $
V^*_{\lambda} \odot V_{\lambda}$ whence this restriction is
non-degenerate, that is, a complex inner product. The operation of
convolution is defined on $ \Bbb C [K^{\Bbb C}]$ and, relative to
the convolution product on $ \Bbb C [K^{\Bbb C}]$,
 the assignment to
$\Phi \in \Bbb C [K^{\Bbb C}]$ of its Fourier coefficient
$\widehat\Phi_{\lambda} \in \roman{End}(V_{\lambda})$ induces a
surjective morphism of algebras
$$
F_{\lambda}\colon\Bbb C [K^{\Bbb C}] \longrightarrow
\roman{End}(V_{\lambda}) \tag2.2
$$
 where
$\roman{End}(V_{\lambda})$ carries its obvious algebra structure,
and this morphism  has the property that, for every $x,y \in
K^{\Bbb C}$ and every $w\in V_{\lambda}$,
$$
(\roman L_x \roman L_y(\Phi))(w)= \TT_{\lambda}(x)
(\Phi((\TT_{\lambda}(y^{-1})w))),\ \Phi \in \Bbb C [K^{\Bbb C}].
\tag2.3
$$
Furthermore, the composite
$$
F_{\lambda} \circ \iota_{\lambda}\colon V_{\lambda}^* \otimes
V_{\lambda} \longrightarrow \roman{End}(V_{\lambda}) \tag2.4
$$
is the canonical isomorphism.

We will use the notation $\alpha = (\alpha_{\lambda}) \in \oplus
\roman{End}(V_{\lambda})$ and $\TT = (\TT_{\lambda}\colon K^{\Bbb
C} \to  \oplus \roman{End}(V_{\lambda}))$, as $\lambda$ ranges
over the highest weights. In terms of this notation, the obvious
action of $K^{\Bbb C}\times K^{\Bbb C}$ on $\oplus
\roman{End}(V_{\lambda})$ is
 given by the association
$$
(x,y,\alpha) \longmapsto \TT(x) \circ \alpha \circ \TT(y^{-1}),\
x,y \in K^{\Bbb C}.
$$
We now recall the algebraic analogue of the Peter-Weyl theorem;
see e.~g. Section 5 of \cite\tspritwo\ for details.

\proclaim{Proposition 2.5}  {\rm (i)} The complex vector space
$\Bbb C [K^{\Bbb C}]$ of representative functions decomposes as
the direct sum
$$
\Bbb C [K^{\Bbb C}] = \oplus_{\lambda}  V^*_{\lambda} \odot
V_{\lambda} \tag2.5.1
$$ of $(K^{\Bbb C}\times K^{\Bbb
C})$-representations and, relative to the complex symmetric
bilinear form {\rm (2.1)},
 the
decomposition is orthogonal.
\newline\noindent
 {\rm (ii)} For each  $\lambda \in \widehat{K^{\Bbb C}}$, the summand
$V^*_{\lambda} \odot V_{\lambda}$ is the isotypical summand of
$\Bbb C [K^{\Bbb C}]$ determined by $\lambda$, and the restriction
of the complex symmetric bilinear form {\rm (2.1)} to this summand
is non-degenerate.
\newline\noindent
{\rm (iii)} Relative to the convolution product on $\Bbb C
[K^{\Bbb
 C}]$,
 the induced morphism
$$
(F_{\lambda})\colon\Bbb C [K^{\Bbb C}] \longrightarrow
\oplus_{\lambda}\roman{End}(V_{\lambda}) \tag2.5.2
$$
of algebras is an isomorphism of $(K^{\Bbb C}\times K^{\Bbb
C})$-representations and yields the decomposition of the
convolution algebra $\Bbb C [K^{\Bbb C}]$ into minimal two-sided
ideals.  \qed
\endproclaim

\beginsection 3. The square integrability of the representative functions

The aim of the present section is to establish the
square-integrability of the representative functions on $K^{\Bbb
C}$ and to reduce the calculation of the requisite integrals over
 $K^{\Bbb C}$ to integrals over $K$.

\proclaim{Lemma 3.1} Each representative function on $K^{\Bbb C}$
is square integrable relative to the measure $\roman e^{-\kappa/t}
\eta \varepsilon$.
\endproclaim

We shall exploit the following integration formula
$$
\int_{\frak k} f(Y) dY = \int_{C^+} \prod_{\alpha \in
R^+}\alpha(Y)^2\left\{\int_{K\big/T} f(\roman{Ad}_y(Y))d(yT)
\right\} dY, \tag3.2
$$
valid for any integrable continuous function $f$ on $\frak k$.
 For the
special case where $\frak k=\frak {su}(2)$, the formula comes
essentially down to integration on $\Bbb R^3$ in ordinary
spherical polar coordinates. For the general case, see e.~g.
\cite\helbotwo\ (Theorem I.5.17, p.~195) or \cite\duiskolk\
((3.14.2) on p.~185 combined with (3.14.4) on p.~187).

Let $\lambda$ be a highest weight. We will use the notation
$\varphi^{\Bbb C}$ etc. for representative functions on $K^{\Bbb
C}$ in the isotypical summand $V^*_{\lambda} \odot V_{\lambda}$ of
$\Bbb C[K^{\Bbb C}]$ associated with $\lambda$ and, accordingly,
we will denote the restriction of $\varphi^{\Bbb C}$ to $K$ by
$\varphi$; then $\varphi$ is necessarily a representative function
on $K$ which lies in the isotypical summand of $L^2(K,dx)$
associated with $\lambda$ by virtue of the ordinary Peter-Weyl
theorem. Lemma 3.1 is implied by the following.

\proclaim{Lemma 3.3} Given the representative function
$\varphi^{\Bbb C}$ on $K^{\Bbb C}$ in the isotypical summand
$V^*_{\lambda} \odot V_{\lambda}$ of $\Bbb C[K^{\Bbb C}]$
associated with the highest weight $\lambda \in \widehat{K^{\Bbb
C}}$,
$$
\int_{K^{\Bbb C}} \overline \varphi^{\Bbb C} \varphi^{\Bbb C}
\roman e^{-\kappa/t} \ETAA \varepsilon = C_{t,\lambda} \int_K
\overline \varphi \varphi dx, \quad C_{t,\lambda}
=(t\pi)^{\dim(K)/2}\roman e^{t|\lambda+\rho|^2}.
$$
\endproclaim

To prepare for the proof, we will denote by $\chi^{\Bbb
C}_{\lambda}$ the holomorphic character of $K^{\Bbb C}$ associated
with the highest weight $\lambda$ and, accordingly, we denote by
$\chi_{\lambda}$ the restriction of $\chi^{\Bbb C}_{\lambda}$ to
$K$; this is plainly the irreducible character of $K$ associated
with $\lambda$.

\proclaim{Lemma 3.4} The character $\chi^{\Bbb C}_{\lambda}$ of
the irreducible representation $\TT_{\lambda}\colon K^{\Bbb C} \to
\roman{End}(V_{\lambda})$ of $K^{\Bbb C}$ associated with the
highest weight $\lambda$ satisfies the identity
$$
\int_{K^{\Bbb C}}||\chi^{\Bbb C}_{\lambda}||^2 \roman
e^{-\kappa}\eta\varepsilon = \frac 1{d_{\lambda}}\int_{K^{\Bbb
C}}|| \TT_\lambda||^2
 \roman e^{-\kappa}\eta\varepsilon
$$
where, as before, $d_{\lambda} = \dim(V_{\lambda})$.
\endproclaim

To prepare for the proof of this Lemma recall that,
 given an
$L^2$-function $f$ on $K$, the appropriate version of the {\it
Plancherel\/} theorem says that a function $f$ on $K$ satisfying
suitable hypotheses, e.~g. \lq f smooth\rq\ suffices, admits the
{\it Fourier\/} decomposition
$$
f(x) = \sum_{\theta} d_{\theta}\roman{tr}(\widehat
f_{\theta}\TT_{\theta}(x)),\ x \in K, \tag3.5
$$
where $\theta$ ranges over the highest weights; see e.~g.
\cite\kirilthr\ (2.3.10). Furthermore, one version of the {\it
Plancherel\/} formula takes the form
$$
\frac 1{\roman{vol}(K)}\int_{K}|f(x)|^2dx
=\sum_{\theta}d_{\theta}||\widehat f_{\theta}||^2; \tag3.6
$$
see e.~g. \cite\kirilthr\ (2.3.11).

For $f=\chi_{\lambda}$, the only non-zero Fourier coefficient
equals $\widehat f_{\lambda}=\frac
1{d_{\lambda}}\roman{Id_{V_{\lambda}}}$, and the Fourier
decomposition of the character $\chi_{\lambda}$ takes the form
$$
\chi_{\lambda}(x) = d_{\lambda}\roman{tr}(\widehat
f_{\lambda}\TT_\lambda (x)),\ x \in K.
$$

\demo{Proof of Lemma {\rm 3.4}} Because the measure $\xi=\roman
e^{-\kappa}\eta\varepsilon$ is $K$-bi-invariant it is in
particular invariant under right translation by elements of $K$.
Hence, for every function $f$ on $K^{\Bbb C}$ which is square
integrable relative to this measure, for each $x \in K$,
$$
\int_{K^{\Bbb C}}||f (y)||^2 d\xi(y) = \int_{K^{\Bbb C}}
||f(yx)||^2 d\xi(y).
$$
Integrating this identity over $K$ yields
$$
\roman{vol}(K) \int_{K^{\Bbb C}} ||f(y)||^2 d\xi(y)= \int_{K^{\Bbb
C}}\int_K ||f (yx)||^2 dx d\xi(y).
$$
Given $y\in K^{\Bbb C}$, the Fourier coefficient $\widehat
f^y_{\lambda}$ of the function $f^y$ on $K$ defined by
$$
f^y(x) = \chi^{\Bbb C}_{\lambda}(yx) = \roman{tr}(\TT_\lambda(y)
\TT_{\lambda}(x)),\ x \in K,
$$
is given by
$$
\widehat f^y_{\lambda} = \TT_{\lambda}(y) \widehat f_{\lambda},
$$
and this is the only non-zero coefficient. Hence, given $y \in
K^{\Bbb C}$, applying the Plancherel formula (3.6) on $K$ to the
function $f^y$, we find
$$
\frac 1{\roman{vol}(K)} \int_K ||f^y(x)||^2 dx =d_{\lambda}||
\TT_{\lambda} \widehat f_{\lambda}||^2 = d_{\lambda}|| \frac
1{d_{\lambda}}\TT_{\lambda}(y) ||^2 = \frac 1{d_{\lambda}}||
\TT_\lambda(y)||^2 .
$$
Consequently
$$
 \int_{K^{\Bbb C}}||\chi^{\Bbb C}_{\lambda}(y)||^2 d\xi(y)=
\frac 1{\roman{vol}(K)}\int_{K^{\Bbb C}}\int_K ||f^y(x)||^2 dx
d\xi(y) = \frac 1{d_{\lambda}}\int_{K^{\Bbb C}}||
\TT_\lambda(y)||^2 d\xi(y)
$$
as asserted. \qed
\enddemo

\demo{Proof of Lemma {\rm 3.3}} We establish the statement of the
Lemma for the special case where $t=1$. The general case is
reduced to the special case by a change of variables.

As a $(K\times K)$-representation, the isotypical summand
$V^*_{\lambda} \odot V_{\lambda}$ is generated by the character
$\chi^{\Bbb C}_{\lambda}$. Hence it suffices to establish the
assertion for $\varphi^{\Bbb C} = \chi^{\Bbb C}_{\lambda}$. By
Lemma 3.4, it suffices to compute the integral $\int_{K^{\Bbb
C}}|| \TT_\lambda||^2 \roman e^{-\kappa}\eta\varepsilon$. To
compute this integral, let $y = x\, \roman{exp}(iY)$ where as
before $x \in K$ and $Y \in \frak k$. Let $\TT_{\lambda}'\colon
\frak k^{\Bbb C} \to \roman {End}(V_{\lambda})$ denote the
corresponding Lie algebra representation and let $A(Y)\in \roman
{End}(V_{\lambda})$ be given by $A(Y)= i \TT_{\lambda}'(Y))$. Then
$$
\TT_{\lambda}(y) = \TT_{\lambda}(x)
\TT_{\lambda}({\roman{exp}(iY)}) = \TT_{\lambda}(x) \roman{exp}(i
\TT_{\lambda}'(Y))=\TT_{\lambda}(x) \roman{e}^{A(Y)}
$$
and
$$
\TT^*_{\lambda}(y) \TT_{\lambda}(y) =
\left(\roman{e}^{A(Y)}\right)^* \left(\roman{e}^{A(Y)}\right) .
$$
Since the endomorphism $\TT_{\lambda}'(Y)$ is skew-hermitian, the
endomorphism $A(Y)$ is hermitian, that is, $A(Y)^* = A(Y)$ whence
$$
\left(\roman{e}^{A(Y)}\right)^* \left(\roman{e}^{A(Y)}\right) =
\roman{e}^{A(Y)^*}\roman{e}^{A(Y)}= \roman{e}^{2A(Y)}
$$
and thence
$$ || \roman{e}^{A(Y)}||^2 =
\roman{tr}(\roman{e}^{2A(Y)}) = \roman{tr}(\roman{e}^{A(2Y)}) =
\roman{tr}(\roman{e}^{i\TT_{\lambda}'(2Y)})=\chi^{\Bbb
C}_{\lambda}(\roman{exp}(2iY)).
$$
View $\lambda+\rho$ as a point of $\frak k^*$ via the orthogonal
decomposition $\frak k = \frak t \oplus \frak q^+$ where $\frak
q^+$ is the orthogonal complement of $\frak t$ in $\frak k$, and
let $\Omega_{\lambda+\rho}$ be the coadjoint orbit generated by
$\lambda+\rho$. Given $Y \in \frak k$, {\it Kirillov\/}'s
character formula, evaluated at the point $\roman{exp}(2iY)$,
yields the identity
$$
\roman{vol}(\Omega_{\rho}) j(2iY) \chi^{\Bbb
C}_{\lambda}(\roman{exp}(2iY)) = \roman{vol}(\Omega_{\rho})
\eta(Y) \chi^{\Bbb C}_{\lambda}(\roman{exp}(2iY)) =
\int_{\Omega_{\lambda+\rho}}\roman e^{-2\vartheta(Y)} d
\sigma(\vartheta),
$$
cf. \cite\kirilboo, \cite\kirilthr. Here  $\vartheta$ refers to
the variable on $\Omega_{\lambda+\rho}$ and $d\sigma$ denotes the
symplectic volume form on $\Omega_{\lambda+\rho}$.
 Using the diffeomorphism from
$K/T$ onto $\Omega_{\lambda+\rho}$ which sends $yT$ ($y \in K$) to
$(\roman{Ad}^*_y)^{-1}(\lambda+\rho)$, we rewrite the integral in
the form
$$
\align \int_{\Omega_{\lambda+\rho}}\roman e^{-\vartheta(2Y)} d
\sigma(\vartheta)
&=\frac{\roman{vol}(\Omega_{\lambda+\rho})}{\roman{vol}(K/T)}
\int_{K\big/T}\roman e^{-(\roman{Ad}^*_y)^{-1}(\lambda+\rho)(2Y)}
d (yT)
\\
&= \frac{d_{\lambda}\roman{vol}(\Omega_{\rho})}{\roman{vol}(K/T)}
\int_{K\big/T}\roman e^{-2(\lambda+\rho)(\roman{Ad}_y(Y))} d (yT).
\endalign
$$
Hence
$$
 \eta(Y)|| \roman{e}^{A(Y)}||^2
 =\eta(Y)\roman{tr}(\roman{e}^{i\TT_{\lambda}'(2Y)})
 =  \frac {d_{\lambda}}{\roman{vol}(K/T)}
\int_{K\big/T} \roman e^{-2(\lambda+\rho)(\roman{Ad}_y(Y))} d (yT)
$$
whence,  in view of the integration formula (3.2),
$$
\align \int_{\frak k}|| \TT_{\lambda}(x\,\roman{exp}(iY))||^2
&\roman e^{-|Y|^2}\eta(Y) dY = \int_{\frak k} ||
\roman{e}^{A(Y)}||^2 \roman e^{-|Y|^2}\eta(Y) dY
\\
&=  d_{\lambda} \int_{C^+} \prod_{\alpha \in R^+}\alpha(Y)^2
 \left\{ \int_{K\big/T}\roman e^{-2(\lambda+\rho)(\roman{Ad}_y(Y))} d (yT)
  \right\}
  \roman e^{-|Y|^2}dY
\\
&=  d_{\lambda} \int_{\frak k}
  \roman e^{-2(\lambda+\rho)(Y)-|Y|^2}dY
= d_{\lambda}\pi^{\dim(K)/2} \roman e^{|\lambda + \rho|^2} .
\endalign
$$
Consequently
$$
\int_{K^{\Bbb C}}|| \chi^{\Bbb C}_{\lambda}||^2 \roman
e^{-\kappa}\eta\varepsilon =\frac {1}{d_{\lambda}} \int_{K^{\Bbb
C}}|| \TT_{\lambda}||^2 \roman e^{-\kappa}\eta\varepsilon =
\pi^{\dim(K)/2} \roman e^{|\lambda + \rho|^2} \roman{vol}(K)
$$
as asserted. In particular, $ C_{1,\lambda}=\pi^{\dim(K)/2} \roman
e^{|\lambda + \rho|^2}$. \qed
\enddemo

\medskip\noindent
{\bf 4. The constituents given by integral forms}
\smallskip\noindent
Left and right translation turn $\Cal H L^2(K^{\Bbb C},\roman
e^{-\kappa/t} \ETAA \varepsilon)$ into a unitary $(K \times
K)$-representation. To establish the statement (i) of the
holomorphic Peter-Weyl theorem,
it remains to show that the
decomposition (2.5.1) of the vector space of representative
functions on $K^{\Bbb C}$ into isotypical summands determines the
decomposition of $\Cal H L^2(K^{\Bbb C},\roman e^{-\kappa/t} \ETAA
\varepsilon)$ into isotypical summands, that is to say:

\proclaim{Proposition 4.1} There is no irreducible $(K \times
K)$-summand in $\Cal H L^2(K^{\Bbb C},\roman e^{-\kappa/t} \ETAA
\varepsilon)$ beyond those which come from the decomposition {\rm
(2.5.1)}.
\endproclaim

We shall establish this fact via a geometric argument which is
guided by the principle that quantization commutes with reduction.
Our argument relies on the familiar complete reducibility of a
continuous unitary representation of a compact Lie group on a
Hilbert space, see e.~g. Theorem III.5.10 on p.~142 of
\cite\broetomd.

We now begin with the preparations for the proof of Proposition
4.1. For intelligibility, we will first recall a few standard
facts. Let $B$ be a Borel subgroup (maximal solvable subgroup) of
$K^{\Bbb C}$ containing $T^{\Bbb C}$.
 Let $\lambda \in \frak t^*=\roman{Hom}(\frak t, \Bbb R)$
be an integral form and let
 $\vartheta_{\lambda}\colon T^{\Bbb C}\to \Bbb C^*$ be the
corresponding algebraic character of the complexification $T^{\Bbb
C}$ of $T$. Thus $\vartheta_{\lambda}$ is given by the formula
$$
\vartheta_{\lambda}(\roman{exp}(w)) = \roman
e^{i\,\lambda(w)},\quad w \in \frak t^{\Bbb C}, \tag4.2
$$
and the derivative of the restriction of $\vartheta_{\lambda}$ to
the maximal torus $T$ coincides with $\lambda$. The corresponding
algebraic character of $B$ is given by the composite of
$\vartheta_{\lambda}\colon T^{\Bbb C}\to \Bbb C^*$ with the
projection from $B$ to $T^{\Bbb C}$, and we denote this character
by $\vartheta_{\lambda}\colon B \to \Bbb C^*$ as well. The $\Bbb
C$-linear subspace
$$
{}_{\lambda}\Bbb C[K^{\Bbb C}] = \{\phi \in \Bbb C[K^{\Bbb C}];
\phi(qy) = \vartheta_{\lambda}(y)^{-1}\phi(q),\ q \in K^{\Bbb C},
y \in B\} \tag4.3.left
$$
of $\Bbb C[K^{\Bbb C}]$ inherits an algebraic $K^{\Bbb C}$-action
in an obvious fashion,  the $K^{\Bbb C}$-action being given by the
assignment to $(x,\phi) \in K^{\Bbb C}\times {}_{\lambda}\Bbb
C[K^{\Bbb C}]$ of $x \phi \in {}_{\lambda}\Bbb C[K^{\Bbb C}]$
where $(x\phi)(q) = \phi(x^{-1}q)$ ($q \in K^{\Bbb C}$). Let
${}_{\lambda}\Bbb C$ denote the 1-dimensional complex vector space
of complex numbers, viewed as a 1-dimensional rational
representation of $B$ via $\vartheta_{\lambda}$, more precisely,
as a 1-dimensional rational {\it left\/} $B$-module. In terms of
this notation,
 ${}_{\lambda}\Bbb C[K^{\Bbb C}]$ is the
$K^{\Bbb C}$-representation $\roman{Ind}_{B}^{K^{\Bbb C}}
{}_{\lambda}\Bbb C$, the rational $K^{\Bbb C}$-representation
 which is induced from ${}_{\lambda}\Bbb C$.
Likewise the $\Bbb C$-linear subspace
$$
\Bbb C[K^{\Bbb C}]_{\lambda} = \{\psi \in \Bbb C[K^{\Bbb C}];
\psi(y q) = \vartheta_{\lambda}(y)^{-1}\psi(q),\ q \in K^{\Bbb C},
y \in B\} \tag4.3.right
$$
of $\Bbb C[K^{\Bbb C}]$ inherits the algebraic $K^{\Bbb C}$-action
given by the assignment to $(x,\phi) \in K^{\Bbb C}\times \Bbb
C[K^{\Bbb C}]_{\lambda} $ of $x \phi \in \Bbb C[K^{\Bbb
C}]_{\lambda} $ where $(x\phi)(q) = \phi(qx)$ ($q \in K^{\Bbb
C}$). With respect to the
 1-dimensional rational {\it
right\/} $B$-module  $\Bbb C_{\lambda}$ which is the vector space
of complex numbers, made into a rational $B$-representation via
$\vartheta_{\lambda}$,
 $\Bbb C[K^{\Bbb C}]_{\lambda} $ amounts to the $K^{\Bbb
C}$-representation $\roman{Ind}_{B}^{K^{\Bbb C}} \Bbb
C_{\lambda}$, the rational $K^{\Bbb C}$-representation which is
induced from $\Bbb C_{\lambda}$. The inversion mapping $x \mapsto
x^{-1}$ on $ K^{\Bbb C}$ which sends $x \in K^{\Bbb C}$ to
$x^{-1}$ induces an isomorphism
$$
{}_{\lambda}\Bbb C[K^{\Bbb C}] \longrightarrow \Bbb C[K^{\Bbb
C}]_{-\lambda}\tag4.4
$$
of $K^{\Bbb C}$-representations.

The choice of dominant Weyl chamber $C^+$ in $\frak t$ determines
a Borel subgroup of $K^{\Bbb C}$ which we denote by $B^+$.
Throughout, highest weights will be understood relative to this
Weyl chamber. Given the integral form $\lambda\in \frak t^*$, we
denote the corresponding algebraic character of $B^+$ by
$\vartheta^+_{\lambda}\colon B^+ \to \Bbb C^*$, the resulting
algebraic $K^{\Bbb C}$-representation (4.3.left) by
${}_{\lambda}\Bbb C[K^{\Bbb C}]^+$, and the representation
(4.3.right) by $\Bbb C[K^{\Bbb C}]_{\lambda}^+$. Let $C^-$ be the
Weyl chamber in $\frak t$ which is opposite to $C^+$, that is, the
Weyl chamber arising from interchanging positive and negative
roots, and let $B^-$ be the corresponding Borel subgroup of
$K^{\Bbb C}$ containing $T^{\Bbb C}$. Given the integral form
$\lambda\in \frak t^*$, let $\vartheta^-_{\lambda}\colon B^- \to
\Bbb C^*$ denote the corresponding algebraic character of $B^-$
which is the composite of $\vartheta_{\lambda}\colon T^{\Bbb C}\to
\Bbb C^*$ with the projection from $B^-$ to $T^{\Bbb C}$, and
denote the resulting algebraic $K^{\Bbb C}$-representation
(4.3.left)
 by ${}_{\lambda}\Bbb C[K^{\Bbb C}]^-$ and the corresponding representation
 (4.3.right) by
$\Bbb C[K^{\Bbb C}]_{\lambda}^-$. By construction, precisely when
$\lambda$ lies in the dominant Weyl chamber $C^+$, the complex
vector spaces ${}_{\lambda}\Bbb C[K^{\Bbb C}]^+$ and $\Bbb
C[K^{\Bbb C}]_{\lambda}^-$ are non-zero and the resulting
representations are irreducible algebraic $K^{\Bbb
C}$-representations. Below, to establish Proposition 4.1, we shall
take $V_{\lambda} ={}_{\lambda}\Bbb C[K^{\Bbb C}]^+$ as $\lambda$
ranges over the highest weights.

\proclaim{Proposition 4.5} Let $\lambda$ be a highest weight
(relative to $C^+$). The assignment to
 $(\psi,\phi)\in \Bbb C[K^{\Bbb C}]^-_{\lambda}
\times{}_{\lambda}\Bbb C[K^{\Bbb C}]^+ $ of $\langle
\psi,\phi\rangle (=(\psi
* \phi)(e))$, cf. {\rm (2.1)}, yields a perfect pairing
$$
\Bbb C[K^{\Bbb C}]^-_{\lambda} \otimes{}_{\lambda}\Bbb C[K^{\Bbb
C}]^+ \longrightarrow \Bbb C \tag4.5.1
$$
which induces  isomorphisms
$$
\aligned \Bbb C[K^{\Bbb C}]^-_{\lambda} \longrightarrow
({}_{\lambda}\Bbb C[K^{\Bbb C}]^+)^* &=\roman{Hom}(
{}_{\lambda}\Bbb C[K^{\Bbb C}]^+,\Bbb C)
\\
{}_{\lambda}\Bbb C[K^{\Bbb C}]^+ \longrightarrow (\Bbb C[K^{\Bbb
C}]^-_{\lambda})^* &=\roman{Hom}( \Bbb C[K^{\Bbb
C}]^-_{\lambda},\Bbb C)
\endaligned
\tag4.5.2
$$
of algebraic $K^{\Bbb C}$-representations. \qed
\endproclaim

As a side remark we note that, in view of the Borel-Weil theorem,
the inclusions of the spaces ${}_{\lambda}\Bbb C[K^{\Bbb C}]^+$
and $\Bbb C[K^{\Bbb C}]_{\lambda}^-$ into the corresponding spaces
of {\it holomorphic\/} functions on $K^{\Bbb C}$ come down to
identity mappings, that is, there is no difference between
algebraic and holomorphic functions at this point.

To recall the familiar descriptions of the spaces
${}_{\lambda}\Bbb C[K^{\Bbb C}]^+$ and $\Bbb C[K^{\Bbb
C}]_{\lambda}^-$ in terms of complex line bundles, let for the
moment $\lambda\in \frak t^*$ be a general integral form. Consider
the complex line bundle
$$
{}_{\lambda}\BETA^{\pm} \colon K^{\Bbb C} \times_{B^{\pm}}
{}_{\lambda}\Bbb C @>>> K^{\Bbb C}\big/B^{\pm} \tag4.6.left
$$
on $K^{\Bbb C}\big/B^{\pm}$. By construction,
 the assignment to $\phi \in {}_{\lambda}\Bbb C[K^{\Bbb C}]^+$
of the induced algebraic section $s_{\phi}$ of
${}_{\lambda}\BETA^{+}$ yields an isomorphism
$$
{}_{\lambda}\Bbb C[K^{\Bbb C}]^+ @>>>
\Gamma_{\roman{alg}}({}_{\lambda}\BETA^{+}) \tag4.7$+$
$$
of complex vector spaces, and in this fashion
$\Gamma_{\roman{alg}}({}_{\lambda}\BETA^+)$ acquires the structure
of an algebraic  $K^{\Bbb C}$-representation. In view of the
Borel-Weil theorem, $\Gamma_{\roman{alg}}({}_{\lambda}\BETA^+)$ is
non-zero when $\lambda$ lies in the dominant Weyl chamber $C^+$.
Likewise, given the general integral form $\lambda$,  consider the
complex line bundle
$$
\BETA^{\pm}_{\lambda} \colon  \Bbb C_{\lambda}\times _{B^{\pm}}
K^{\Bbb C} @>>> B^{\pm}\backslash K^{\Bbb C} \tag4.6.right
$$
on $B^{\pm}\backslash K^{\Bbb C}$. By construction, the assignment
to $\psi \in \Bbb C[K^{\Bbb C}]_{\lambda}^-$ of the induced
algebraic section $s_{\psi}$ of $\BETA^-_{\lambda}$ yields an
isomorphism
$$
\Bbb C[K^{\Bbb C}]_{\lambda}^- @>>>
\Gamma_{\roman{alg}}(\BETA^-_{\lambda}) \tag4.7$-$
$$
of complex vector spaces and, in this fashion, the space
$\Gamma_{\roman{alg}}(\BETA^-_{\lambda})$ of algebraic sections of
$\BETA^-_{\lambda}$ acquires the structure of an algebraic
$K^{\Bbb C}$-representation. For a general integral form
$\lambda$, the algebraic mapping
$$
K^{\Bbb C} \times_{B^\pm} {}_{\lambda}\Bbb C @>>> \Bbb
C_{-\lambda} \times _{B^\pm} K^{\Bbb C},\quad (x,v) \longmapsto
(v,x^{-1}),
$$
where $x \in K^{\Bbb C}$ and $v \in \Bbb C$, induces an
isomorphism
$$
{}_{-\lambda}\BETA^{\pm}  @>>> \BETA^{\pm}_{\lambda} \tag4.8
$$
of algebraic line bundles   which, on the bases, is the algebraic
isomorphism
$$
K^{\Bbb C}\big /B^\pm \longrightarrow  B^\pm \backslash K^{\Bbb
C},\quad xB^\pm \mapsto B^\pm x^{-1},\  x \in K^{\Bbb C},
$$
and this isomorphism induces an isomorphism of algebraic $K^{\Bbb
C}$-representations
$$
\Gamma_{\roman{alg}}({}_{-\lambda}\BETA^\pm)
\longrightarrow\Gamma_{\roman{alg}}(\BETA^\pm_{\lambda}) \tag4.9
$$
which is plainly compatible with the isomorphism (4.4) of
algebraic $K^{\Bbb C}$-representations between ${}_{-\lambda}\Bbb
C[K^{\Bbb C}]^\pm$ and $\Bbb C[K^{\Bbb C}]^\pm_{\lambda}$. In
particular, by the Borel-Weil theorem,
$\Gamma_{\roman{alg}}(\BETA^-_{\lambda})$ is non-zero precisely
when $-\lambda$ lies in the Weyl chamber corresponding to $B^-$,
that is, when $\lambda$ lies in the dominant Weyl chamber
determined by $B^+$
 and, in this case,
$\Bbb C[K^{\Bbb
C}]_{\lambda}^-\cong\Gamma_{\roman{alg}}(\BETA^-_{\lambda})$ is an
irreducible $K^{\Bbb C}$-representation with highest weight
$w(-\lambda +\rho)-\rho$ where $w$ is the unique element of the
Weyl group such that $w(-\lambda +\rho)-\rho$ lies in the interior
of the dominant Weyl chamber $C^+$ where as before $\rho$ refers
to one half the sum of the positive roots.

For $\nu \in \frak k^*$, let  $\Cal O_{\nu} = K\nu\subseteq \frak
k^*$ be the coadjoint orbit generated by $\nu$. For completeness,
we recall that the $K$-action on $\frak k^*$ is given in the
standard way, that is, by means of the association
$$
K\times \frak k^* \longrightarrow \frak k^*, \quad (x,\chi)
\mapsto x\chi = \roman{Ad}^*_{x^{-1}}(\chi),\quad x \in K, \, \chi
\in \frak k^* .
$$
Let $\lambda$ be an integral form in $C^+$. Then $\lambda+\rho$
lies in the interior of the dominant Weyl chamber, and the
 coadjoint orbit $\Cal O_{\lambda+\rho}$ of $\lambda+\rho$
in $\frak k^*$ has maximal dimension, whether or not the orbit of
$\lambda$ has maximal dimension, that is, the stabilizer of the
point  ${\lambda+\rho}$ of $\frak k^*$ is minimal and coincides
with the maximal torus $T$; likewise the stabilizer of the point
${-(\lambda+\rho)}$ of $\frak k^*$ is minimal and coincides with
the maximal torus $T$. Since the inclusion mapping $K \subseteq
K^{\Bbb C}$ induces a diffeomorphism $K\big/T \to K^{\Bbb C}\big
/B^{+}$, the assignment to $x \in K$ of
$$
x(\lambda+\rho)=\roman{Ad}^*_{x^{-1}} (\lambda + \rho) \in \frak
k^* \tag4.10$+$
$$
induces an embedding
$$
{}_{\lambda}\mu^+ \colon K^{\Bbb C}\big /B^+ \longrightarrow \frak
k^* \tag4.11$+$
$$
of the homogeneous space $K^{\Bbb C}\big /B^+$ into $\frak k^*$
which induces a $K$-equivariant diffeomorphism from $K^{\Bbb
C}\big /B^+$ onto the coadjoint orbit $\Cal O_{\lambda+\rho}$. It
is well known that the Kirillov-Kostant-Souriau symplectic
structure $\sigma_{\lambda+\rho}$ on $\Cal O_{\lambda+\rho}$
combines with the complex structure on $K^{\Bbb C}\big /B^+$ to a
positive $K$-invariant K\"ahler structure on both $\Cal
O_{\lambda+\rho}$ and $K^{\Bbb C}\big /B^+$ in such a way that
${}_{\lambda}\mu^+$ identifies the two resulting K\"ahler
manifolds in a $K$-equivariant fashion and such that
${}_{\lambda}\mu^+$ is a $K$-equivariant momentum mapping.
Furthermore, relative to the additional structure on $K^{\Bbb
C}\big /B^+$, the line bundle ${}_{\lambda}\beta^+$ is positive,
in fact, a prequantum bundle, by construction necessarily
$K$-equivariant, the unique hermitian connection being the
requisite connection. Likewise the inclusion mapping $K \subseteq
K^{\Bbb C}$ induces a diffeomorphism $T\backslash K \to
B^-\backslash K^{\Bbb C}$, and the assignment to $y \in K$ of
$$
-y^{-1}(\lambda+\rho)=-\roman{Ad}^*_{y} (\lambda + \rho) \in \frak
k^* \tag4.10$-$
$$
induces an embedding
$$
\mu^-_{\lambda} \colon B^-\backslash K^{\Bbb C} \longrightarrow
\frak k^* \tag4.11$-$
$$
of the homogeneous space $B^-\backslash K^{\Bbb C}$ into $\frak
k^*$. For convenience, we convert the obvious $K^{\Bbb C}$-action
on the right of $B^-\backslash K^{\Bbb C}$ in the standard way to
a left action via the association
$$
K^{\Bbb C} \times (B^-\backslash K^{\Bbb C}) \longrightarrow
B^-\backslash K^{\Bbb C},\ (y,B^-x) \mapsto B^- x y^{-1},\ x,y \in
K^{\Bbb C}. \tag4.12
$$
With this convention, (4.11$-$) induces a $K$-equivariant
diffeomorphism from $B^-\backslash K^{\Bbb C}$ onto the coadjoint
orbit $\Cal O_{-(\lambda+\rho)}$ in such a way that (i) the
Kirillov-Kostant-Souriau symplectic structure
$\sigma_{-(\lambda+\rho)}$ on $\Cal O_{-(\lambda+\rho)}$ combines
with the complex structure on $B^-\backslash K^{\Bbb C}$ to a
positive $K$-invariant K\"ahler structure on both $\Cal
O_{-(\lambda+\rho)}$ and $B^-\backslash K^{\Bbb C}$,  that (ii)
$\mu^-_{\lambda}$ identifies the two resulting K\"ahler manifolds
in a $K$-equivariant fashion,  that (iii) $\mu^-_{\lambda}$ is a
$K$-equivariant momentum mapping, and such that (iv) relative to
the additional structure on $B^-\backslash K^{\Bbb C}$, the line
bundle $\beta^-_{\lambda}$ is positive, in fact, a prequantum
bundle, by construction necessarily $K$-equivariant, the requisite
connection being the unique hermitian connection.

As before, we consider $\roman T^*K$ as a Hamiltonian $(K\times
K)$-space relative to the $(K\times K)$-action which arises from
the lifts of the left translation and of the right translation
action on $K$. The momentum mapping
$$
\mu^{K \times K} \colon \roman T^*K @>>> \frak k^* \times \frak
k^*
$$
for this $(K\times K)$-action on $\roman T^*K \cong K^{\Bbb C}$ is
well known to be given by  the association
$$
\roman T^*K \ni \alpha_x \longmapsto (\alpha_x \circ \roman
R_x,\alpha_x \circ \roman L_x) \in \frak k^* \times \frak k^*,\ x
\in K, \, \alpha_x \in \roman T_x^*K ,
$$
where $\roman R_x \colon \frak k =\roman T_eK \to \roman T_xK$ and
$\roman L_x \colon \frak k =\roman T_eK \to \roman T_xK$ refer to
the operations of left- and right translation, respectively, by
$x\in K$. With an abuse of notation, we denote the corresponding
momentum mapping on $K^{\Bbb C}$ by
$$
\mu^{K \times K} \colon K^{\Bbb C} @>>> \frak k^* \times \frak k^*
$$
as well, and we denote the symplectic structure on $K^{\Bbb C}$ by
$\sigma_K$.

Consider the product manifold
$$
N^{\times}= K^{\Bbb C} \times (K^{\Bbb C}\big/B^+) \times
(B^-\backslash K^{\Bbb C}),
$$
endowed with the product K\"ahler structure. Let $\sigma^{\times}$
be the resulting product symplectic structure which underlies the
product K\"ahler structure, essentially the sum of $\sigma_K$,
$\sigma_{\lambda+\rho}$, and  $\sigma_{-(\lambda+\rho)}$. The
group $K^{\Bbb C}\times K^{\Bbb C}$ acts on $N^{\times}$ in the
obvious fashion, that is, the action on $K^{\Bbb C}$ is given by
left- and right translation, that on $K^{\Bbb C}\big/B^+$ by the
projection to the first factor $K^{\Bbb C}$, followed by the
$K^{\Bbb C}$-action on $K^{\Bbb C}\big/B^+$, and that on
$B^-\backslash K^{\Bbb C}$ by the projection to the second factor
$K^{\Bbb C}$, followed be the $K^{\Bbb C}$-action on
$B^-\backslash K^{\Bbb C}$. Furthermore, by construction, the
symplectic structure $\sigma^{\times}$ is $(K \times
K)$-invariant. Let
$$
\mu^{\times}\colon  N^{\times}@ >>> \frak k^* \times \frak k^*
$$
be the $(K\times K)$-momentum mapping for the $(K\times K)$-action
on $N^{\times}$ relative to the symplectic structure
$\sigma^{\times}$. This momentum mapping is essentially the sum of
the momentum mapping $\mu^{ K \times K}$ and the two momentum
mappings (4.11$+$) and (4.11$-$). The $(K\times K)$-reduced space
$\left(\mu^{\times}\right)^{-1}(0,0)\big/(K\times K)$ at the point
zero of $\frak k^* \times \frak k^*$ boils down to a single point.

We will denote the complex vector space  of holomorphic functions
on $K^{\Bbb C}$ by $\Cal H(K^{\Bbb C})$. Left and right
translation on $K^{\Bbb C}$ turn $\Cal H (K^{\Bbb C})$ into a
holomorphic $(K^{\Bbb C}\times K^{\Bbb C})$-representation. By
construction, the product line bundle
$$
 \BETA^{\times}=\BETA
\times {}_{\lambda}\BETA^+ \times \BETA^-_{\lambda}
$$
is a holomorphic $(K^{\Bbb C}\times K^{\Bbb C})$-equivariant line
bundle and, in view of the isomorphisms (4.7$+$) and (4.7$-$),
since the complex vector spaces ${}_{\lambda}\Bbb C[K^{\Bbb C}]^+$
and $\Bbb C[K^{\Bbb C}]_{\lambda}^-$ are finite-dimensional, as a
$(K^{\Bbb C}\times K^{\Bbb C})$-representation, the space of
holomorphic sections of this line bundle amounts to the tensor
product
$$
{}_{\lambda}\Bbb C[K^{\Bbb C}]^+ \otimes \Bbb C[K^{\Bbb
C}]_{\lambda}^- \otimes \Gamma_{\roman{hol}} (\BETA)\cong
{}_{\lambda}\Bbb C[K^{\Bbb C}]^+ \otimes \Bbb C[K^{\Bbb
C}]_{\lambda}^- \otimes \Cal H(K^{\Bbb C})
$$
of representations. In view of the isomorphisms (4.5.2), as a
$(K^{\Bbb C}\times K^{\Bbb C})$-representation, this tensor
product may be written as
$$
\roman{Hom}_{\Bbb C}\left( \Bbb C[K^{\Bbb C}]_{\lambda}^- \otimes
{}_{\lambda}\Bbb C[K^{\Bbb C}]^+,\Cal H(K^{\Bbb C})\right).
$$
By construction, the product line bundle $\BETA^{\times}$ is a
holomorphic $(K\times K)$-equivariant prequantum bundle on the
K\"ahler manifold $N^{\times}$. The $(K\times K)$-reduced space
$\left(\mu^{\times}\right)^{-1}(0,0)\big/(K\times K)$ at the point
zero of $\frak k^* \times \frak k^*$ is a single point. Indeed,
consider the point $(e,B^+, B^-)$ of $N^{\times}$. This point lies
in $\left(\mu^{\times}\right)^{-1}(0,0)$, and the $(K\times
K)$-orbit of this point is the entire zero locus
$\left(\mu^{\times}\right)^{-1}(0,0)$. This observation implies at
once that the $(K\times K)$-reduced space
$\left(\mu^{\times}\right)^{-1}(0,0)\big/(K\times K)$ at the point
zero of $\frak k^* \times \frak k^*$ is a single point. For
completeness we note that the stabilizer of the point $(e,B^+,
B^-)$ of $\left(\mu^{\times}\right)^{-1}(0,0)$ is a copy of the
maximal torus $T$. With these preparations out of the way, we
conclude that the space
 of $(K\times K)$-invariant
holomorphic sections of the product line bundle $\BETA^{\times}$
is at most 1-dimensional, that is, the space
$$
\roman{Hom}_{\Bbb C}\left( \Bbb C[K^{\Bbb C}]_{\lambda}^- \otimes
{}_{\lambda}\Bbb C[K^{\Bbb C}]^+,\Cal H(K^{\Bbb
C})\right)^{K\times K} \tag4.13
$$
is at most 1-dimensional.

We explain briefly under somewhat more general circumstances how
one arrives at the last conclusion: Let $G$ be a compact Lie
group, let $N$ be a $G$-Hamiltonian K\"ahler manifold, with
$G$-equivariant momentum mapping $\mu \colon N \to \frak g^*$, and
suppose that $G$ preserves the complex structure on $N$. Then $G$
preserves the associated Riemannian metric as well, and the
$G$-action extends canonically to a holomorphic $G^{\Bbb
C}$-action on $N$. The {\it saturation\/} of the zero locus
$\mu^{-1}(0)$ is the subspace $G^{\Bbb C}\mu^{-1}(0) \subseteq N$,
and the inclusion $\mu^{-1}(0)\subseteq G^{\Bbb C}\mu^{-1}(0)$
induces a homeomorphism from the reduced space $N_0 =\mu^{-1}(0)
\big/ G$ onto the $G^{\Bbb C}$-quotient $G^{\Bbb C}\mu^{-1}(0)
\big /G^{\Bbb C}$. In this fashion, the reduced space $N_0$
acquires a complex analytic structure.
 Let $\zeta \colon E \to N$ be a
$G$-invariant prequantum bundle, and let
$$
\zeta^0 \colon E|_{G^{\Bbb C}\mu^{-1}(0)} \longrightarrow G^{\Bbb
C}\mu^{-1}(0)
$$
be the restriction of $\zeta$ to $G^{\Bbb C}\mu^{-1}(0) \subseteq
N$. Passing to $G^{\Bbb C}$-quotients, we obtain the coherent
analytic  sheaf $ \zeta_0 \colon E_0 \longrightarrow N_0 $ on
$N_0$, not necessarily an ordinary line bundle. The canonical
morphism $ \pi \colon \Gamma(\zeta^0)^G \longrightarrow
\Gamma(\zeta_0) $ of complex vector spaces is plainly injective
(even an isomorphism, but we do not need this fact): A
$G$-equivariant section of $\zeta^0$ inducing the zero section of
$\zeta_0$ is manifestly the zero section. Furthermore, the
restriction mapping
 from $\Gamma(\zeta)^G$ onto $\Gamma(\zeta^0)^G$ is an isomorphism.
Applying this reasoning to $N=N^{\times}$ and $G=K\times K$, we
conclude that the vector space (4.13) is at most 1-dimensional as
asserted. For intelligibility we note that, as a space, the
saturation
$$
(K^{\Bbb C}\times K^{\Bbb C})\left(\mu^{\times}\right)^{-1}(0,0)
\subseteq N^{\times}
$$
amounts to a homogeneous space of the kind $(K^{\Bbb C}\times
K^{\Bbb C})\big/T^{\Bbb C}$, the complexification $T^{\Bbb C}$ of
the maximal torus $T$ of $K$ being suitably embedded into $K^{\Bbb
C}\times K^{\Bbb C}$, but we shall not need this fact. However,
this observation shows that, for topological reasons, the
saturation cannot be all of $N^{\times}$.

We now take $V_{\lambda} ={}_{\lambda}\Bbb C[K^{\Bbb C}]^+$. Then,
in view of Proposition 4.5, $V^*_{\lambda} \cong \Bbb C[K^{\Bbb
C}]{}_{\lambda}^-$. Since we already know that, in the
decomposition (2.5.1) of the vector space $\Bbb C[K^{\Bbb C}]$ of
representative functions on $K^{\Bbb C}$, $V^*_{\lambda}\odot
V_{\lambda}$ is the isotypical summand corresponding to $\lambda$,
and since, by virtue of Lemma 3.1, $V^*_{\lambda}\odot
V_{\lambda}$ is actually a subspace of the Hilbert space $\Cal H
L^2(K^{\Bbb C},\roman e^{-\kappa/t} \ETAA \varepsilon)$, we
conclude that the space
$$
\roman{Hom}_{\Bbb C}\left( \Bbb C[K^{\Bbb C}]_{\lambda}^- \otimes
{}_{\lambda}\Bbb C[K^{\Bbb C}]^+, \Cal H L^2(K^{\Bbb C},\roman
e^{-\kappa/t} \ETAA \varepsilon)\right)^{K\times K}
$$
is 1-dimensional. However, this vector space is that of morphisms
of $(K\times K)$-representations from $\Bbb C[K^{\Bbb
C}]_{\lambda}^- \otimes {}_{\lambda}\Bbb C[K^{\Bbb C}]^+$ to $\Cal
H L^2(K^{\Bbb C},\roman e^{-\kappa/t} \ETAA \varepsilon)$. Since
this space is 1-dimensional, it is generated by a single such
morphism, and this morphism picks out the $(K\times
K)$-irreducible constituent $\Bbb C[K^{\Bbb C}]_{\lambda}^-
\otimes {}_{\lambda}\Bbb C[K^{\Bbb C}]^+\cong \roman{End}_{\Bbb
C}({}_{\lambda}\Bbb C[K^{\Bbb C}]^+)$ from $\Cal H L^2(K^{\Bbb
C},\roman e^{-\kappa/t} \ETAA \varepsilon)$. In other words,
$V^*_{\lambda}\odot V_{\lambda}$ is the isotypical summand in
$\Cal H L^2(K^{\Bbb C},\roman e^{-\kappa/t} \ETAA \varepsilon)$
determined by $\lambda$.

These observations imply that  the vector space $\Bbb C[K^{\Bbb
C}]$ of representative functions on $K^{\Bbb C}$ is dense in the
Hilbert space $ \Cal H L^2(K^{\Bbb C},\roman e^{-\kappa/t} \ETAA
\varepsilon)$. This proves Proposition 4.1 and hence establishes
statement (i) of the holomorphic Peter-Weyl theorem, Theorem 1.14.

\medskip\noindent {\bf 5. The abstract identification with the vertically
polarized Hilbert space}
\smallskip\noindent
The vertically polarized Hilbert space arising from geometric
quantization on $\roman T^*K$ is a Hilbert space of half forms.
Haar measure $dx$ on $K$ then yields a concrete realization of
this Hilbert space as $L^2(K,dx)$. In this section we will compare
the Hilbert space $ \Cal H L^2(K^{\Bbb C},\roman e^{-\kappa/t}
\ETAA \varepsilon)$ with the vertically polarized Hilbert space.
This will in particular provide a proof of statement (ii) of the
holomorphic Peter-Weyl theorem.

Let $\lambda$ be a highest weight. Let $W_{\lambda}$ denote the
space of complex representative functions on $K$ which arise by
restriction to $K$ of the holomorphic functions in $V_{\lambda}$.
Since a holomorphic function on $K^{\Bbb C}$ is determined by its
values on $K$, this restriction mapping is the identity mapping of
complex vector spaces, in fact, of $K$-representations. To justify
the distinction in notation,  we note that the embedding
$\iota_{\lambda}$ given as (1.8) above yields an embedding
$$
\iota_{\lambda} \colon W^*_{\lambda} \otimes W_{\lambda}
\longrightarrow R(K)=\Bbb C[K^{\Bbb C}]
$$
and, maintaining the notation $\odot$ introduced in Section 2, we
write
$$
W^*_{\lambda} \odot W_{\lambda}= \iota_{\lambda}( W^*_{\lambda}
\otimes W_{\lambda}) \subseteq R(K).
$$
The $K$-representation $W^*_{\lambda} \otimes W_{\lambda}$
inherits a $K$-invariant inner product from the embedding into
$L^2(K,dx)$. On the other hand, $V^*_{\lambda} \otimes
V_{\lambda}$ acquires an inner product from its embedding  into
the Hilbert space $\Cal H L^2(K^{\Bbb C},\roman e^{-\kappa/t}
\ETAA \varepsilon)$ induced by (1.8) which turns $V^*_{\lambda}
\otimes V_{\lambda}$ into a unitary $K$-representation, but the
relationship between the inner products on $V^*_{\lambda} \otimes
V_{\lambda}$ and $W^*_{\lambda} \otimes W_{\lambda}$ is not a
priori clear. We therefore distinguish the resulting unitary
$K$-representations $W_{\lambda}$ and $V_{\lambda}$ in notation as
indicated.

Let $\langle \,\cdot\, , \, \cdot \,\rangle_K$ denote the
normalized inner product on $L^2(K,dx)$ given by
$$
\langle f,h\rangle_K= \frac 1{\roman{vol}(K)} \int_K\overline f h
dx. \tag5.1
$$
As usual, we endow ${\oplus} _{\lambda \in \widehat{K^{\Bbb
C}}}\roman{End}(W_\lambda)$ with the  inner product which, on the
summand $\roman{End}(W_\lambda)$, is given by
$$
 {d_{\lambda}}\langle \,\cdot\, , \,\cdot\,
\rangle_{\lambda}.
$$
This inner product differs from the inner product (1.13); see the
completion of the proof of the holomorphic Peter-Weyl theorem
given below for an explanation. Then $\widehat{\oplus} _{\lambda
\in \widehat{K^{\Bbb C}}}\roman{End}(W_\lambda)$ refers to the
completion relative to this inner product. As in the situation of
the inner product (1.13), up to a constant, the resulting norm on
each $\roman{End}(W_\lambda)$ coincides with the {\it
Hilbert-Schmidt\/} norm. For ease of exposition, we spell out the
ordinary {\it Peter-Weyl\/} theorem in the following form.

\proclaim{Proposition 5.2} {\rm (i)} The space $R(K)$ of
representative functions on $K$ is dense in $L^2(K,dx)$ and, as a
unitary $(K\times K)$-representation, $L^2(K,dx)$ decomposes as
the direct sum
$$
L^2(K,dx) = \widehat{\oplus} _{\lambda} (W^*_{\lambda} \odot
W_{\lambda})
 \cong \widehat{\oplus}
_{\lambda}\roman{End}(W_\lambda) \tag5.2.1
$$
into $(K\times K)$-isotypical summands as $\lambda$ ranges over
the highest weights.
\newline\noindent
{\rm (ii)} Relative to the convolution product $*$ on $L^2(K,dx)$,
as $\lambda$ ranges over the highest weights, the assignment to an
$L^2$-function $f$ on $K$ of its Fourier coefficients $\widehat
f_{\lambda}\in \roman{End}(W_\lambda)$ yields an isomorphism
$$
L^2(K,dx) \longrightarrow \widehat{\oplus}
_{\lambda}\roman{End}(W_\lambda) \tag5.2.2
$$
of Hilbert algebras where $L^2(K,dx)$ is endowed with the
normalized inner product $\langle \,\cdot\, , \, \cdot
\,\rangle_K$.
\endproclaim

The following is an immediate consequence of the ordinary and the
holomorphic Peter-Weyl theorem, combined with the explicit
determination of the constants $C_{t,\lambda}$ for the highest
weights $\lambda$ given in Lemma 3.3, viz.
$C_{t,\lambda}=(t\pi)^{\dim(K)/2}\roman e^{t|\lambda+\rho|^2}$.

\proclaim{Theorem 5.3} The association
$$
V^*_{\lambda} \odot V_{\lambda} \ni \varphi^{\Bbb C} \longmapsto
C_{t,\lambda}^{1/2}\varphi= (t\pi)^{\dim(K)/4}\roman
e^{t|\lambda+\rho|^2/2} \varphi\in W^*_{\lambda} \odot
W_{\lambda},
$$
as $\lambda$ ranges over the highest weights, induces a unitary
isomorphism
$$
H_t\colon\Cal H L^2(K^{\Bbb C},\roman e^{-\kappa/t} \ETAA
\varepsilon) \longrightarrow  L^2(K,dx)  \tag5.3.1
$$
of unitary $(K\times K)$-representations. \qed
\endproclaim

\demo{Completion of the proof  of the holomorphic Peter-Weyl
theorem} Let
$$
H^{\roman{End}}_t \colon \widehat{\oplus} _{\lambda \in
\widehat{K^{\Bbb C}}}\roman{End}(V_\lambda)@>>> \widehat{\oplus}
_{\lambda \in \widehat{K^{\Bbb C}}}\roman{End}(W_\lambda)
$$
be the obvious unitary isomorphism of $(K\times
K)$-representations which, restricted to the summand
$\roman{End}(V_\lambda)$, is given by multiplication by
$C_{t,\lambda}^{1/2}$, as $\lambda$ ranges over the highest
weights. By construction, the diagram
$$
\CD
 \Cal H L^2(K^{\Bbb C},\roman e^{-\kappa/t} \ETAA
\varepsilon) @>{H_t}>>L^2(K,dx)
\\
@VVV @VVV
\\
\widehat{\oplus} _{\lambda \in \widehat{K^{\Bbb
C}}}\roman{End}(V_\lambda)@>>{H^{\roman{End}}_t}> \widehat{\oplus}
_{\lambda \in \widehat{K^{\Bbb C}}}\roman{End}(W_\lambda)
\endCD
$$
is commutative where the unlabelled vertical arrows are given by
the assignment to a function of its Fourier coefficients.
Moreover, in view of Theorem 5.3, the upper horizontal arrow is an
isomorphism of unitary $(K\times K)$-representations, the lower
horizontal arrow is such an isomorphism as just pointed out
 and,
by virtue of the ordinary Peter-Weyl theorem, the right-hand
vertical arrow is an isomorphism of Hilbert algebras. In view of
the algebraic version of the Peter-Weyl theorem, Proposition 2.5
above, we conclude that the convolution product on the algebra
$\Bbb C[K^{\Bbb C}]$ extends to a convolution product on $\Cal H
L^2(K^{\Bbb C},\roman e^{-\kappa/t} \ETAA \varepsilon)$ and that
the left-hand vertical arrow is an isomorphism of Hilbert algebras
as asserted.~\qed
\enddemo

As a consequence of the holomorphic Peter-Weyl theorem, we will
now spell out a holomorphic version of the {\it Plancherel\/}
theorem. Given the holomorphic function $\Phi$ on $K^{\Bbb C}$, we
refer to the series $ \sum_{\lambda} d_{\lambda}\roman{tr}
\left(\TT_{\lambda}(y)\widehat \Phi_{\lambda}\right) $ in the
variable $y\in K^{\Bbb C}$ as the {\it holomorphic Fourier
series\/} of $\Phi$. Up to a change of variable, the holomorphic
Fourier series of $\Phi$ coincides with the ordinary Fourier
series of the restriction of $\Phi$ to $K$. We will denote by
$||\,\cdot\,||_{t,K^{\Bbb C}}$ the norm associated with the inner
product (1.5).

\proclaim {Corollary 5.4}{\rm [Holomorphic Plancherel theorem]}
The holomorphic Fourier series of a holomorphic function $\Phi$ on
$K^{\Bbb C}$ that is square integrable relative to the measure
$\roman e^{-\kappa/t}\eta \varepsilon$ converges to $\Phi$ in
$\Cal HL^2(K^{\Bbb C},\roman e^{-\kappa/t}\eta \varepsilon)$ and
hence converges to $\Phi$ pointwise as well. Furthermore, given
the family $\left(c_{\lambda}\right)_{\lambda \in \widehat
{K^{\Bbb C}}}$ where $c_{\lambda} \in \roman{End}(V_{\lambda})$,
the series $ \sum_{\lambda} d_{\lambda}\roman{tr}
\left(\TT_{\lambda}(y)c_{\lambda}\right) $ furnishes a holomorphic
function on $K^{\Bbb C}$ which is square-integrable relative to
the measure $\roman e^{-\kappa/t}\eta \varepsilon$
 if and only if
the series $\sum d_{\lambda}C_{t,\lambda}||c_{\lambda}||^2$
converges; if this happens to be the case, when $\Phi$ denotes the
resulting holomorphic function, the Plancherel formula takes the
form
$$
||\Phi||^ 2_{t,K^{\Bbb C}}=\frac 1{\roman{vol}(K)}\int_{K^{\Bbb
C}}|\Phi|^2 \roman e^{-\kappa} \eta\varepsilon = \sum
d_{\lambda}C_{t,\lambda}||c_{\lambda}||^2. \tag5.4.1
$$
\endproclaim

\demo{Proof} Let $\lambda$ be a highest weight, let
$\TT_{\lambda}\colon K^{\Bbb C} \to \roman{End}(V_{\lambda})$ be
the associated irreducible rational representation of $K^{\Bbb
C}$, and let
$$
\TT^{\Bbb
C}_{\lambda}=\frac{\TT_{\lambda}}{C_{t,\lambda}^{1/2}},\quad
\widehat \Phi^{\Bbb C}_{\lambda}=C_{t,\lambda}^{1/2} \widehat
\Phi_{\lambda}.
$$
Then, with the obvious extension of the notation $\langle\,\cdot\,
, \,\cdot\, \rangle_{t,K^{\Bbb C}}$, we have
$$
 \widehat \Phi^{\Bbb C}_{\lambda} = \left\langle
{\TT^{\Bbb C} _{\lambda}},\Phi\right\rangle_{t,K^{\Bbb C}} = \frac
1{\roman{vol}(K)} \int_{K^{\Bbb C}}{\overline \TT^{\Bbb C}
_{\lambda}}\Phi \roman e^{-\kappa} \eta\varepsilon,
$$
that is, $\widehat \Phi^{\Bbb C}_{\lambda}$ is the Fourier
coefficient of $\Phi$ relative to $\lambda$, calculated with
respect to the Hilbert space $\Cal HL^2(K^{\Bbb C},\roman
e^{-\kappa/t}\eta \varepsilon)$. Thus the holomorphic Fourier
series of $\Phi$ can be written in the form
$$
\sum_{\lambda} d_{\lambda}\roman{tr} \left({\TT^{\Bbb C}
_{\lambda}(y)}\widehat \Phi^{\Bbb C}_{\lambda}\right)
\left(=\sum_{\lambda} d_{\lambda}\roman{tr} \left(\frac
{\TT_{\lambda}(y)}{C_{t,\lambda}^{1/2}}C_{t,\lambda}^{1/2}\widehat
\Phi_{\lambda}\right)\right) .
$$
Given the representative function $\varphi^{\Bbb C}$ in
$V^*_{\lambda} \otimes V_{\lambda}$, under the isomorphism
(5.3.1), the representative function ${\varphi^{\Bbb
C}}/{C_{t,\lambda}^{1/2}}$ on $K^{\Bbb C}$ goes to the restriction
$\varphi$ of $\varphi^{\Bbb C}$ to $K$. These observations imply
the assertions. \qed
\enddemo

\noindent{\smc Remark 5.5.} A version of the holomorphic
Plancherel Theorem may be found in Lemmata 9 and 10 of
\cite\bhallthr. According to \cite\cerezone\ (Proposition 12), the
holomorphic Fourier series (referred to in \cite\cerezone\ as a
Fourier-Laurent series) of a general holomorphic function  on
$K^{\Bbb C}$, not necessarily square integrable relative to the
measure $\roman e^{-\kappa/t}\eta \varepsilon$, converges
uniformly on compact sets. This fact has been extended to
holomorphic functions on the complexification of a general
symmetric space of a compact Lie group in \cite\lassatwo\ (Theorem
3 in Subsection 5.5). The statement of Corollary 5.4 can, perhaps,
be deduced from the estimates given in \cite\lassaone\ but we do
not know whether this has been worked out in the literature.
Corollary 5.4 includes the convergence in the Hilbert space $\Cal
HL^2(K^{\Bbb C},\roman e^{-\kappa/t}\eta \varepsilon)$; for the
convergence in this Hilbert space,  see also Lemma 10 in
\cite\bhallthr. In Theorem 9 (iii) of \cite\bhallthr\ a formula
similar to (5.4.1) above is given, valid relative to any
sufficiently regular bi-invariant measure.

\medskip\noindent {\bf 6. The Blattner-Kostant-Sternberg pairing}
\smallskip\noindent
In this section we will show that the isomorphism
(5.3.1) is realized by the B(lattner)K(ostant)S(ternberg)-pairing,
multiplied by a global constant; see e.~g. \cite\sniabook,
\cite\woodhous\ for details on the BKS-pairing. We maintain the
notation $\langle \,\cdot\, , \, \cdot \,\rangle_{t,K^{\Bbb C}}$
for the normalized inner product on $\Cal H L^2(K^{\Bbb C},\roman
e^{-\kappa/t} \ETAA \varepsilon)$ induced by the measure $\roman
e^{-\kappa/t} \ETAA \varepsilon$.

 Let $\Phi$ be a holomorphic function on $K^{\Bbb
C}$ which is square integrable relative to $\roman
e^{-\kappa/t}\ETAA \varepsilon$ and let $\FFF$ be a square
integrable function on $K$; the ordinary BKS-pairing
$\langle\,\cdot\, , \,\cdot \, \rangle_{\roman{BKS}}$ between the
two half-form Hilbert spaces $\Cal H L^2(K^{\Bbb C},\roman
e^{-\kappa/t} \ETAA \varepsilon)$ and $L^2(K,dx)$ assigns the
integral
$$
\langle \Phi,\FFF\rangle_{\roman{BKS}} = \frac
1{\roman{vol}(K)}\int_K \int_{\frak k} \overline{\Phi(x\,
\roman{exp}(iY))} \FFF(x) \roman e^{- \frac {|Y|^2}{2t}}
\eta(Y/2)dY dx \tag6.1.1
$$
to $\Phi$ and $\FFF$ provided this integral exists. The requisite
calculation which yields the explicit form (6.1.1) of the
BKS-pairing under the present circumstances is given in the
appendix of \cite\bhallone, where the notation $ \zeta(Y) =
\eta(Y/2) $ is used ($Y\in \frak k$).

We will now show that (6.1.1) extends to a pairing which is
defined everywhere, that is, to a pairing of the kind
$$
\langle\,\cdot\, , \,\cdot \, \rangle_{\roman{BKS}} \colon \Cal H
L^2(K^{\Bbb C},\roman e^{-\kappa/t}\ETAA \varepsilon) \otimes
L^2(K,dx) \longrightarrow \Bbb C. \tag6.1.2
$$
We do not assert that the integral is absolutely convergent for
every $\Phi$ and $F$, though. To begin with we note that is
manifest that, given a holomorphic function $\Phi$ on $K^{\Bbb C}$
which is square integrable relative to $\roman e^{-\kappa/t}\ETAA
\varepsilon$, when the complex function $F_{\Phi}$ on $K$ given by
the expression
$$
(F_\Phi)(x) = \int_{\frak k} {\Phi(x \,\roman {exp}{(iY)})} \roman
e^{- \frac {|Y|^2}{2t}} \eta(Y/2)  dY, \quad x \in K, \tag6.2
$$
is well defined, that is, when the integral exists for every $x \in K$,
$$
\langle \Phi,\FFF\rangle_{\roman{BKS}} = \langle
F_\Phi,\FFF\rangle_{K}. \tag6.3
$$

\proclaim{Lemma 6.4} Let $\lambda$ be a highest weight, let
$\varphi^{\Bbb C}$ be a representative function on $K^{\Bbb C}$ in
the isotypical summand $V^*_{\lambda} \odot V_{\lambda}$ of $\Cal
H L^2(K^{\Bbb C},\roman e^{-\kappa/t} \ETAA \varepsilon)$
associated with $\lambda$ and, as before, let $\varphi$ denote the
restriction of $\varphi^{\Bbb C}$ to $K$, necessarily a
representative function on $K$ which lies in the isotypical
summand $W^*_{\lambda} \odot W_{\lambda}$ of $L^2(K,dx)$. Then the
integral {\rm (6.2)} exists for every $x\in K$, and the resulting
function $F_{\varphi^{\Bbb C}}$ on $K$ is given by
$$
F_{\varphi^{\Bbb C}} = D_{t,\lambda} \varphi, \quad D_{t,\lambda}
=(2t\pi)^{\dim(K)/2}\roman e^{t{|\lambda+\rho|^2}/2}.
$$
\endproclaim

\demo{Proof} We establish the statement of the Lemma for the
special case where $t=1$. The general case is reduced to the
special case by a change of variables.

As a $(K\times K)$-representation, the isotypical summand
$V^*_{\lambda} \odot V_{\lambda}$ associated with $\lambda$ is
spanned by the character $\chi^{\Bbb C}_{\lambda}$ of $K^{\Bbb C}$
associated with the highest weight $\lambda$. Thus it suffices to
establish the claim for $\varphi^{\Bbb C}=\chi^{\Bbb
C}_{\lambda}$, and we will now do so:

In view of the integration formula (3.2), given $x \in K$,
$$
\align F_{\chi^{\Bbb C}_{\lambda}}(x) &= \int_{\frak k}
{\chi^{\Bbb C}_{\lambda}(x \,\roman
{exp}{(iY)})} \roman e^{- {|Y|^2}/2} \eta(Y/2)  dY \\
 &= \frac 1{\roman{vol}(T)}\int_{C^+}
\prod_{\alpha \in R^+}\alpha(Y)^2\left\{\int_{K} \chi^{\Bbb
C}_{\lambda}(x\, \roman{exp}(\roman{Ad}_y(iY)))dy \right\} \roman
e^{- {|Y|^2}/2} \eta(Y/2) dY
\\
&= \frac 1{\roman{vol}(T)}\int_{C^+} \prod_{\alpha \in
R^+}\alpha(Y)^2\left\{\int_K \chi^{\Bbb C}_{\lambda}(y^{-1}xy\,
\roman{exp}(iY))dy \right\}\roman e^{- {|Y|^2}/2} \eta(Y/2) dY.
\endalign
$$
Let $x \in K$ and $Y \in \frak k$; using the formula
$$
\int_K \chi^{\Bbb C}_{\lambda}(y^{-1}xy\, \roman{exp}(iY))dy
=\frac{\roman{vol}(K)}{d_{\lambda}}\chi_{\lambda}(x) \chi^{\Bbb
C}_{\lambda}(\roman{exp}(iY))
$$
where, as before, $d_{\lambda}$ denotes the dimension of the
irreducible representation  associated with $\lambda$, we conclude
$$
F_{\chi^{\Bbb C}_{\lambda}}(x)=
 \frac {\roman{vol}(K\big/T)}{d_{\lambda}}
{\chi_{\lambda}(x)} \int_{C^+} \prod_{\alpha \in
R^+}\alpha(Y)^2\chi^{\Bbb C}_{\lambda}( \roman{exp}(iY))\roman
e^{- {|Y|^2}/2} \eta(Y/2)dY .
$$
Given $Y \in \frak k$, {\it Kirillov\/}'s character formula, cf.
\cite\kirilboo, \cite\kirilthr, evaluated at the point
$\roman{exp}(iY)$, yields the identity
$$
\roman{vol}(\Omega_{\rho}) j(iY) \chi^{\Bbb
C}_{\lambda}(\roman{exp}(iY)) =
 \roman{vol}(\Omega_{\rho})\eta(Y/2)
 \chi^{\Bbb C}_{\lambda}(\roman{exp}(iY)) =
\int_{\Omega_{\lambda+\rho}}\roman e^{-\vartheta(Y)} d
\sigma(\vartheta).
$$
Now, as in the proof of Lemma 3.3, using the diffeomorphism from
$K/T$ onto $\Omega_{\lambda+\rho}$ which sends $yT$ ($y \in K$) to
$(\roman{Ad}^*_y)^{-1}(\lambda+\rho)$, we rewrite the integral as
an integral over $K\big/T$ and obtain the identity
$$
 \eta(Y/2)
\chi^{\Bbb C}_{\lambda}(\roman{exp}(iY)) =
\frac{d_{\lambda}}{\roman{vol}(K/T)} \int_{K\big/T} \roman
e^{-(\lambda+\rho)(\roman{Ad}_y(Y))} d (yT).
$$
Hence
$$
\align F_{\chi^{\Bbb C}_{\lambda}}(x) &=
 \frac {\roman{vol}(K\big/T)}{d_{\lambda}}
{\chi_{\lambda}(x)} \int_{C^+} \prod_{\alpha \in
R^+}\alpha(Y)^2\chi^{\Bbb C}_{\lambda}( \roman{exp}(iY))\roman
e^{- {|Y|^2}/2} \eta(Y/2)dY
\\
&= {\chi_{\lambda}(x)}
 \int_{C^+} \prod_{\alpha
\in R^+}\alpha(Y)^2 \left\{ \int_{K\big/T}\roman
e^{-(\lambda+\rho)(\roman{Ad}_y(Y))} d (yT)\right\} \roman e^{-
{|Y|^2}/2} dY
\\
&= {\chi_{\lambda}(x)} \int_{\frak k}\roman e^{-(\lambda+\rho)(Y)-
{|Y|^2}/2} dY = (2 \pi)^{\dim(K)/2} \roman
e^{{|\lambda+\rho|^2}/2} {\chi_{\lambda}(x)}
\endalign
$$
whence, in particular, $ D_{1,\lambda}= (2 \pi)^{\dim(K)/2} \roman
e^{{|\lambda+\rho|^2}/2}$ as asserted. \qed
\enddemo

\proclaim{Theorem 6.5} The {\rm BKS}-pairing {\rm (6.1.1)} extends
to a (non-degenerate) $(K\times K)$-invariant pairing of the kind
{\rm (6.1.2)}. Furthermore, the assignment to a representative
function $\Phi$ on $K^{\Bbb C}$ of the function $F_\Phi$ on $K$
induces a bounded $(K\times K)$-equivariant operator
$$
\Theta_t \colon \Cal H L^2(K^{\Bbb C},\roman e^{-\kappa/t}\ETAA
\varepsilon) \longrightarrow L^2(K,dx) \tag6.5.1
$$
such that
$$
\langle \Phi,\FFF\rangle_{\roman{BKS}} = \langle
\Theta_t(\Phi),\FFF\rangle_{K} \tag6.5.2
$$
and such that, when $\varphi^{\Bbb C}$ is a member of the
isotypical summand $V^*_{\lambda} \odot V_{\lambda}$,
$$
\Theta_t(\varphi^{\Bbb C})= F_{\varphi^{\Bbb C}} = D_{t,\lambda}
\varphi, \tag6.5.3
$$
where as before $\varphi$ refers to the restriction of
$\varphi^{\Bbb C}$ to $K$. Finally, the operator
$$
(4t\pi)^{-\dim(K)/4} \Theta_t \colon \Cal H L^2(K^{\Bbb C},\roman
e^{-\kappa/t} \ETAA \varepsilon) \longrightarrow L^2(K,dx)
\tag6.5.4
$$
sends a representative function $\varphi^{\Bbb C} \in
V^*_{\lambda} \odot V_{\lambda}$ to
$C_{t,\lambda}^{1/2}\varphi=(t\pi)^{\dim(K)/4}\roman
e^{t|\lambda+\rho|^2/2}\varphi$ and thus coincides with the
unitary isomorphism {\rm(5.3.1)} of $(K\times K)$-representations.
\endproclaim

\demo{Proof} Let $\lambda_1$ and $\lambda_2$ be two highest
weights, let $\varphi^{\Bbb C}$ be a representative function on
$K^{\Bbb C}$ which is a member of the isotypical summand
$V^*_{\lambda_1} \odot V_{\lambda_1}$ associated with the highest
weight $\lambda_1$, and let $\psi$ be a representative function on
$K$ which is a member of the isotypical summand $W^*_{\lambda_2}
\odot W_{\lambda_2}$ associated with the highest weight
$\lambda_2$. In view of the identity (6.3) and Lemma 6.4,
$$
\langle \varphi^{\Bbb C},\psi\rangle_{\roman{BKS}} = \langle
F_{\varphi^{\Bbb C}},\psi\rangle_{K} = D_{t,\lambda_1} \langle
\varphi,\psi\rangle_{K}.
$$
Hence, by virtue of the ordinary Peter-Weyl theorem and of the
holomorphic Peter-Weyl theorem, the BKS-pairing (6.1.2) is
everywhere defined. By construction,  the pairing is
$K$-bi-invariant.

Let $\varphi^{\Bbb C}$ be a representative function on $K^{\Bbb
C}$ which is a member of the isotypical summand $V^*_{\lambda}
\odot V_{\lambda}$ associated with the highest weight $\lambda$.
Since
$$
\frac {C_{t,\lambda}}{D_{t,\lambda}^2}
=\frac{(t\pi)^{\dim(K)/2}\roman e^{t|\lambda+\rho|^2}}
{{(2t\pi)^{\dim(K)}}\roman
e^{t|\lambda+\rho|^2}}=((4t\pi)^{-\dim(K)/4})^2 ,
$$
by virtue of Lemma 3.3 and Lemma 6.4,
$$
\int_{K^{\Bbb C}} \overline \varphi^{\Bbb C} \varphi^{\Bbb C}
\roman e^{-\kappa/t} \ETAA \varepsilon = \frac
{C_{t,\lambda}}{D_{t,\lambda}^2} \int_K \overline
{F_{\varphi^{\Bbb C}}} F_{\varphi^{\Bbb C}} dx =
((4t\pi)^{-\dim(K)/4})^2 \int_K \overline {F_{\varphi^{\Bbb C}}}
F_{\varphi^{\Bbb C}} dx.
$$
In view of the ordinary Peter-Weyl theorem and of the holomorphic
Peter-Weyl theorem, this identity implies the remaining assertions
of Theorem 6.5. \qed \enddemo

Let $\Theta_t^* \colon L^2(K,dx) \to \Cal H L^2(K^{\Bbb C},\roman
e^{-\kappa/t} \ETAA \varepsilon)$ be the adjoint of $\Theta_t$.
Let $\lambda$ be a highest weight, let $\varphi^{\Bbb C}\in
V_{\lambda}^*\odot V_{\lambda}$ and let $\varphi\in
W_{\lambda}^*\odot W_{\lambda}$ be the restriction of
$\varphi^{\Bbb C}$ to $K$. Define the number $A_{t,\lambda}$ by
$\Theta_t^*(\varphi) = A_{t,\lambda}\varphi^{\Bbb C}$. Then
$$
D_{t,\lambda}\langle \varphi, \varphi\rangle_K =\langle
\Theta_t(\varphi^{\Bbb C}), \varphi\rangle_K = \langle
\varphi^{\Bbb C}, \Theta_t^*\varphi\rangle_{t,K^{\Bbb C}}
 = A_{t,\lambda} \langle \varphi^{\Bbb C},
 \varphi^{\Bbb C}\rangle_{t,K^{\Bbb C}} =
A_{t,\lambda}C_{t,\lambda} \langle \varphi, \varphi\rangle_K
$$
whence $ A_{t,\lambda} = \frac  {D_{t,\lambda}}{C_{t,\lambda}} =
2^{\dim(K)/2} \roman e^{-t{|\lambda+\rho|^2}/2}$. Hence
$$
\Theta_t^*(\varphi) = 2^{\dim(K)/2} \roman
e^{-t{|\lambda+\rho|^2}/2}\varphi^{\Bbb C}. \tag6.6
$$

\proclaim{Corollary 6.7} The resulting operator
$$
(4t\pi)^{-\dim(K)/4} \Theta_t^* \colon L^2(K,dx) \longrightarrow
\Cal H L^2(K^{\Bbb C},\roman e^{-\kappa/t} \ETAA \varepsilon)
\tag6.7.1
$$
is unitary and coincides with the inverse of the isomorphism {\rm
(5.3.1)}.
\endproclaim

\noindent{\smc Remark 6.8.} As explained already in the
introduction, the unitarity of the BKS-pairing map, multiplied by
a suitable constant, has been established in \cite\bhallone\ by
means of the heat kernel techniques developed in \cite\bhallthr.
In \cite\liuwanhu, using the very same heat kernel methods, the
authors have shown that the unitarity of the BKS-pairing map can
be reduced to a computation on matrix entries. Likewise, the proof
of Theorem 5.3 reduces the abstract unitary equivalence between
the two Hilbert spaces involved to inspection of certain square
integrals of representative functions. However, the proof of
Theorem 5.3 is direct and independent of heat kernel techniques,
and in fact the statement of Theorem 5.3 is at first independent
of the BKS-pairing map as well.

\beginsection 7. The spectral decomposition of the energy operator

Let $\Delta_K$ denote the {\it Casimir\/} operator on $K$
associated with the bi-invariant Riemannian metric on $K$. When
$X_1,\dots,X_m$ is an orthonormal basis of $\frak k$,
$$
\Delta_K = X^2_1 +\dots + X^2_m
$$
in the universal algebra $\roman U(\frak k)$ of $\frak k$, cf.
e.~g. \cite\enelson\ (p.~591). The Casimir operator depends only
on the Riemannian metric, though. Since the metric on $K$ is
bi-invariant, so is the operator $\Delta_K$; hence, by Schur's
lemma, each isotypical summand $W^*_{\lambda}\odot W_\lambda
\subseteq L^2(K,dx)$ is an eigenspace, whence the representative
functions are eigenfunctions for $\Delta_K$. The eigenvalue of
$\Delta_K$ corresponding to the highest weight $\lambda$ is known
to be given explicitly by $-\varepsilon_{\lambda}$ where
$$
\varepsilon=(|\lambda+\rho|^2-|\rho|^2),
$$
cf. e.~g.
\cite\helbotwo\ (Ch. V.1 (16) p.~502). The present sign is
dictated by the interpretation in terms of the energy given below.
Thus $\Delta_K$ acts on each isotypical summand
$W^*_{\lambda}\odot W_\lambda$ as scalar multiplication by
$-\varepsilon_{\lambda}$. The Casimir operator is known to
coincide with the nonpositive Laplace-Beltrami operator associated
with the (bi-invariant) Riemannian metric on $K$, see e.~g.
\cite\taylothr\ (A 1.2). In the Schr\"odinger picture (vertical
quantization on $\roman T^*K$), the unique extension $\widehat
E_K$ of the operator $-\frac 12 \Delta_K$  to an unbounded
self-adjoint operator on $L^2(K,dx)$ is the quantum mechanical
{\it energy\/} operator associated with the Riemannian metric,
whence the spectral decomposition of this operator refines in the
standard manner to the Peter-Weyl decomposition of $L^2(K,dx)$
into isotypical $(K\times K)$-summands.

Via the embedding of $\frak k$ into $\frak k^{\Bbb C}$, the
operator $\Delta_K$ is a differential operator on $K^{\Bbb C}$. In
view of the holomorphic Peter-Weyl theorem, the unitary transform
(5.3.1) (or, equivalently, (6.7.1),) is compatible with the
operator $\Delta_K$. Consequently, in the holomorphic quantization
on $\roman T^*K\cong K^{\Bbb C}$, via the transform (5.3.1) (or,
equivalently, via the BKS-pairing map (6.5.1) multiplied by
$(4t\pi)^{-\dim(K)/4}$), the operator $\widehat E_{K^{\Bbb C}}$
which arises as the unique extension of the operator $-\frac 12
\Delta_K$ on $\Cal H L^2(K^{\Bbb C},\roman e^{-\kappa/t} \ETAA
\varepsilon)$ to an unbounded self-adjoint operator is the quantum
mechanical {\it energy\/} operator associated with the Riemannian
metric, and the spectral decomposition of this operator refines to
the holomorphic Peter-Weyl decomposition of $\Cal H L^2(K^{\Bbb
C},\roman e^{-\kappa/t} \ETAA \varepsilon)$ into isotypical
$(K\times K)$-summands.

Finally, we note that, in terms of the Casimir operator
$\Delta_K$, the identity (6.6) may plainly be written in the form
$$
\Theta_t^*(\varphi) = 2^{\dim(K)/2} \roman e^{-t|\rho|^2/2} \roman
e^{-t{\Delta_K}/2} \varphi^{\Bbb C}, \tag7.1
$$
where $\varphi$ is any representative function on $K$.
 In this description of the operator $\Theta_t^*$,
the highest weights, present in the description (6.6), no longer
appear explicitly. Consequently, for any smooth function $f$ on
$K$, $\Theta_t^*(f)$ is the unique holomorphic function on
$K^{\Bbb C}$ whose restriction to $K$ is given by
$$
\Theta_t^*(f)|_K =  2^{\dim(K)/2} \roman e^{-t|\rho|^2/2} \roman
e^{t{\Delta_K}/2}f. \tag7.2
$$
In Theorem 2.6(1) of \cite\bhallone, this operator $\Theta_t^*$ is
written as $\Pi_{\hbar}$, where the parameter $\hbar$ corresponds
to the present notation $t$. The value $\roman e^{t{\Delta_K}/2}f$
is also given by
$$
(\roman e^{t{\Delta_K}/2}f)(y) =\int_K p_t (yx^{-1}) f(x) dx =
(p_t*f)(y), \ y\in K, \tag7.3
$$
where $p_t$ is the fundamental solution of the heat equation
$\frac {du}{dt} = \frac 12 \Delta_K(u)$ on $K$, subject to the
initial condition that $p_0$ be the Dirac distribution supported
at $e\in K$ \cite\bhallone, \cite\enelson \ (Section 8).

\smallskip
\noindent{\smc Remark 7.4.} With some computational effort, the
numerical values of the eigenvalues $-\varepsilon_{\lambda}$ of
the Laplace operator being known, the abstract isomorphism between
the two Hilbert spaces spelled out in Theorem 5.3 above can also
be derived from Theorem 10 in \cite\bhallthr\ which, in turn, is
proved via heat kernel techniques. Needless, perhaps, to repeat
again: Our approach to the abstract isomorphism between the two
Hilbert spaces spelled out in Theorem 5.3 is independent of heat
kernel techniques.

\medskip\noindent {\bf 8. Relationship with the naive Hilbert
space} \smallskip\noindent We refer to the Hilbert space $\Cal H
L^2(K^{\Bbb C},\roman e^{-\kappa /t}\varepsilon)$ of holomorphic
functions that are square-integrable relative to the measure
$\roman e^{-\kappa /t}\varepsilon$ as the {\it naive\/} Hilbert
space. We content ourselves with the following simplified version
of the corresponding holomorphic Peter-Weyl theorem.

\proclaim{Proposition 8.1} The Hilbert space $\Cal H L^2(K^{\Bbb
C},\roman e^{-\kappa /t}\varepsilon)$ contains the vector space
which underlies the algebra $\Bbb C[K^{\Bbb C}]$ of representative
functions on $K^{\Bbb C}$ as a dense subspace.
\endproclaim

\demo{Proof} Since the measure is Gaussian, standard arguments
involving the appropriate estimates show that each representative
function is square-integrable relative to the measure $\roman
e^{-\kappa /t}\varepsilon$. The reasoning which establishes
Proposition 4.1 is also valid for the naive Hilbert space. This
completes the proof. \qed
\enddemo

For the highest weight $\lambda$, define the constant $\widetilde
C_{t,\lambda}$ by the identity
$$
\int_{K^{\Bbb C}}||\chi^{\Bbb C}_{\lambda}||^2 \roman
e^{-\kappa}\varepsilon = \widetilde C_{t,\lambda} \roman{vol}(K).
$$
Analogously to Theorem 5.3, we now have the following.

\proclaim{Theorem 8.2} The association
$$
V^*_{\lambda} \odot V_{\lambda} \ni \varphi^{\Bbb C} \longmapsto
\widetilde C_{t,\lambda}^{1/2}\varphi\in W^*_{\lambda} \odot
W_{\lambda},
$$
as $\lambda$ ranges over the highest weights, induces a unitary
isomorphism
$$
\widetilde H_t\colon\Cal H L^2(K^{\Bbb C},\roman e^{-\kappa/t}
\varepsilon) \longrightarrow  L^2(K,dx)  \tag8.2.1
$$
of unitary $(K\times K)$-representations. \qed
\endproclaim

Consequently the two Hilbert spaces $\Cal H L^2(K^{\Bbb C},\roman
e^{-\kappa /t}\eta \varepsilon)$ and $\Cal H L^2(K^{\Bbb C},\roman
e^{-\kappa /t}\varepsilon)$ are unitarily equivalent as $(K\times
K)$-representations. However we do not know how to compute the
values of the constants $\widetilde C_{t,\lambda}$. A tool like
Kirillov's character formula does not seem to be available for
this case. Furthermore, we do not know whether there is a
candidate for a pairing inducing the equivalence between the two
Hilbert spaces.

\medskip

\centerline{\smc References} \smallskip
\widestnumber\key{999}

\ref \no \rbieltwo \by R. Bielawski \paper K\"ahler metrics on
$G^{\Bbb C}$ \jour J. reine angew. Mathematik \vol 559 \yr 2003
\pages 123--136 \finalinfo{\tt math.DG/0202255}
\endref

\ref \no \broetomd \by  Th. Br\"ocker and T. tom Dieck \book
Representations of Compact Lie groups \bookinfo Graduate Texts in
Mathematics, No. 98 \publ Springer Verlag \publaddr Berlin $\cdot$
Heidelberg  $\cdot$   New York $\cdot$  Tokyo \yr 1985
\endref

\ref \no \cerezone \by A. Cerezo \paper Solutions analytiques des
\'equations  invariantes sur un groupe compact ou complexe
r\'eductif \jour Ann. Inst. Fourier \vol 25 \yr 1975 \pages
249--277
\endref

\ref \no \duflotwo \by M. Duflo \paper Op\'erateurs
diff\'erentiels bi-invariants sur un groupe de Lie \jour Ann. Sci.
\'Ecole Norm. Sup. (4) \vol 10 \yr 1977 \pages 265--283
\endref

\ref \no \duiskolk \by J.~J. Duistermaat and J.~A.~C. Kolk \book
{Lie groups}\bookinfo {Universitext} \publ {Springer-Verlag}
\publaddr{Berlin $\cdot$ Heidelberg $\cdot$ New York $\cdot$ Tokyo
}\yr{2000}\endref

\ref \no \farauone \by J. Faraut \paper Espaces hilbertiens
invariants de fonctions holomorphes \paperinfo in: Analyse sur les
groupes de Lie et th\'eorie des repr\'esentations \jour
S\'eminaires et Congr\`es, Soc. Math. France \vol 7 \yr 2003
\pages 101--167
\endref

\ref \no \bhallthr \by  B. C. Hall \paper The Segal-Bargmann
\lq\lq Coherent State\rq\rq\ Transform for compact Lie groups
 \jour Journal of Functional Analysis \vol 122 \yr 1994 \pages
103--151
\endref

\ref \no \bhallfou \by  B. C. Hall \paper
Phase space bounds for quantum mechanics on a compact Lie group
 \jour Comm. in Math. Physics \vol 184 \yr 1997 \pages
233--250
\endref

\ref \no \bhallone \by  B. C. Hall \paper Geometric quantization
and the generalized Segal-Bargmann transform for Lie groups of
compact type \jour Comm. in Math. Physics \vol 226 \yr 2002 \pages
233--268 \finalinfo {\tt quant.ph/0012015}
\endref

\ref \no \helbotwo \by S. Helgason \book Groups and Geometric
Analysis. Integral geometry, invariant differential operators, and
spherical functions \bookinfo Pure and Applied Mathematics, vol.
113 \publ Academic Press Inc. \publaddr Orlando, Fl. \yr 1984
\endref

\ref \no \poiscoho \by J. Huebschmann \paper Poisson cohomology
and quantization \jour J. reine angew.
 Mathematik
\vol  408 \yr 1990 \pages 57--113
\endref

\ref \no  \souriau \by J. Huebschmann \paper On the quantization
of Poisson algebras \paperinfo Symplectic Geometry and
Mathematical Physics, Actes du colloque en l'honneur de Jean-Marie
Souriau, P. Donato, C. Duval, J. Elhadad, G.M. Tuynman, eds.;
Progress in Mathematics, Vol. 99 \publ Birkh\"auser Verlag
\publaddr Boston $\cdot$ Basel $\cdot$ Berlin \yr 1991 \pages
204--233
\endref

\ref \no \kaehler \by J. Huebschmann \paper K\"ahler spaces,
nilpotent orbits, and singular reduction \linebreak
\jour Memoirs
AMS \vol 172/814 \yr 2004 \publ Amer. Math. Society \publaddr
Providence, R. I.\finalinfo {\tt math.DG/0104213}
\endref

\ref \no \lradq \by J. Huebschmann \paper Lie-Rinehart algebras,
descent, and quantization \paperinfo Galois theory, Hopf algebras,
and semiabelian categories \jour Fields Institute Communications
\vol 43 \publ Amer. Math. Society \publaddr Providence, R.~I.
\pages 295--316 \yr 2004 \finalinfo {\tt math.SG/0303016}
\endref

\ref \no \qr \by J. Huebschmann \paper K\"ahler quantization and
reduction \finalinfo {\tt math.SG/0207166} \jour J.  reine angew.
Mathematik
\vol 591
\yr 2006
\pages 75--109
\endref

\ref \no \varna \by J. Huebschmann \paper Classical phase space
singularities and quantization \paperinfo in: Quantum Theory and
Symmetries. IV. V. Dobrev, ed.  \publ Heron Press \publaddr Sofia
\yr 2006 \pages 51--65 \finalinfo {\tt math-ph/0610047}
\endref

\ref \no \bedlepro \by J. Huebschmann \paper Singular
Poisson-K\"ahler geometry of certain adjoint quotients\paperinfo
in: The Mathematical Legacy of C. Ehresmann, J. Kubarski, and R.
Wolak, eds. \jour Banach Center Publications (to appear)
\finalinfo {\tt math.SG/0610614}
\endref

\ref \no \adjoint \by J. Huebschmann \paper Stratified K\"ahler
structures on adjoint quotients \linebreak \finalinfo{\tt
math.DG/0404141}
\endref

\ref \no \hurusch \by J. Huebschmann, G. Rudolph, and M. Schmidt
\paper A gauge model for quantum mechanics on a stratified space
\finalinfo{\tt hep-th/0702017}
\endref

 \ref \no \kirilboo \by A. A. Kirillov
\book Elements of the theory of representations. Translated from
the Russian by E. Hewitt \bookinfo Grundlehren der Mathematischen
Wissenschaften, vol. 220 \publ Springer Verlag \publaddr Berlin
$\cdot$ Heidelberg $\cdot$ New York \yr 1976
\endref

\ref \no \kirilthr \by A. A. Kirillov \paper Merits and demerits
of the orbit method \jour Bull. Amer. Math. Soc. \vol 36 \yr 1999
\pages 433--488
\endref

\ref \no \lassaone \by M. Lassalle \paper Sur la transformation de
Fourier-Laurent dans un groupe analytique complexe r\'eductif
\jour Ann. Institut Fourier \vol 28 \yr 1978 \pages 115--138
\endref

\ref \no \lassatwo \by M. Lassalle \paper S\'eries de Laurent des
fonctions holomorphes dans la complexification d'un espace
sym\'etrique compact \jour Ann. Scient. \'Ecole Normale
Sup\'erieure \vol 11 \yr 1978 \pages 167--210
\endref

\ref \no \lempszoe \by L. Lempert and R. Sz\"oke \paper Global
solutions of the homogeneous complex Monge-Am\-p\`ere equations
and complex structures on the tangent bundle of Riemannian
manifolds \jour Math. Ann. \vol 290 \yr 1991 \pages 689--712
\endref

\ref \no \liuwanhu \by Liu, Wei-ping, Wang, Zheng-dong, and Hu,
Da-peng \paper Segal-Bargmann-Hall transform and geometric
quantization \jour Adv. Math. (China) \vol 32 \yr 2003 \pages
509--511
\endref

\ref \no \enelson \by E. Nelson \paper Analytic vectors \jour Ann.
of Mathematics \vol 70 \yr 1959 \pages  572--615
\endref

\ref \no \sniabook \by J. \'Sniatycki \book Geometric quantization
and quantum mechanics \bookinfo Applied Mathematical Sciences
 No.~30
\publ Springer-Verlag \publaddr Berlin $\cdot$ Heidelberg $\cdot$
New York \yr 1980
\endref

\ref \no \tspritwo \by T.~A. Springer \paper Aktionen reduktiver
Gruppen \paperinfo Algebraische Transformationsgruppen und
Invariantentheorie, DMV Seminar Band 13, H.~P. Kraft, P. Slodowy,
T.~A. Springer, eds. \publ Birkh\"auser Verlag \publaddr Boston
$\cdot$ Basel $\cdot$ Berlin \yr 1989 \pages 4--39
\endref

\ref \no \stenzone \by M.~B. Stenzel \paper The Segal-Bargmann
transform on a symmetric space of compact type \jour Journal of
Functional Analysis \vol 165 \yr 1999 \pages 44--58
\endref

\ref \no \szoekone \by R. Sz\"oke \paper Complex structures on
tangent bundles of Riemannian manifolds \jour Math. Ann. \vol 291
\yr 1991 \pages 409--428
\endref

\ref \no \taylothr \by J.  Taylor \paper The Iwasawa decomposition
and limiting behaviour of Brownian motion on symmetric spaces of
non-compact type \jour Cont. Math. \vol 73 \yr 1988 \pages
303--331
\endref

\ref \no \woodhous \by N. M. J. Woodhouse \book Geometric
quantization \bookinfo Oxford Mathematical Monographs, Second
edition \publ Clarendon Press \publaddr Oxford \yr 1991
\endref

\enddocument